\documentclass[12pt]{article}
\usepackage{graphicx}
\usepackage{amsmath}
\usepackage{amssymb}
\usepackage{theorem}

\numberwithin{equation}{section}

\newtheorem{thm}{Theorem}[section]
\newtheorem{lemma}[thm]{Lemma}
\newtheorem{prop}[thm]{Proposition}
\newtheorem{cor}[thm]{Corollary}
{\theorembodyfont{\rmfamily}
\newtheorem{defn}[thm]{Definition}
\newtheorem{example}[thm]{Example}
\newtheorem{rmk}[thm]{Remark}
}

\newcommand{\qed}{\hfill \mbox{\raggedright \rule{.07in}{.1in}}}
 
\newenvironment{proof}{\vspace{1ex}\noindent{\bf
Proof}\hspace{0.5em}}{\hfill\qed\vspace{1ex}}
\newenvironment{pfof}[1]{\vspace{1ex}\noindent{\bf Proof of
#1}\hspace{0.5em}}{\hfill\qed\vspace{1ex}}

\newcommand{\R}{{\mathbb R}}
\newcommand{\C}{{\mathbb C}}
\newcommand{\Z}{{\mathbb Z}}
\newcommand{\N}{{\mathbb N}}

\newcommand{\dist}{\operatorname{dist}}

\newcommand{\SMALL}{\textstyle}

\title{Decay of Correlations for Slowly Mixing Flows}

\author{Ian Melbourne~\thanks{Department of Mathematics and Statistics,
University of Surrey, Guildford GU2 7XH, UK, ism@math.uh.edu
\newline 2000 {\em Mathematics Subject Classification.} 37A25, 37D25, 37D50}}

\date{6 November 2006}

\begin{document}

\maketitle

\begin{abstract}
We show that polynomial decay of correlations is prevalent
for a class of nonuniformly hyperbolic flows.  These flows are the
continuous time analogue of a class of nonuniformly hyperbolic 
diffeomorphisms for which Young proved polynomial decay of correlations.
Roughly speaking, in situations where the decay rate $O(1/n^{\beta})$ has 
previously been proved for diffeomorphisms, we establish the decay rate 
$O(1/t^\beta)$ for typical flows.
Applications include certain classes of semidispersing billiards, as well as
dispersing billiards with vanishing curvature.

In addition, we obtain results for suspension flows with unbounded roof 
functions.   
This includes the planar periodic Lorentz flow with infinite horizon.
\end{abstract}

\section{Introduction}

Dolgopyat~\cite{Dolgopyat98b} has shown that uniformly hyperbolic (Axiom~A)
flows typically mix rapidly, faster than any polynomial rate, for
sufficiently smooth observables.
The restriction to typical flows is necessary; 
there exist uniformly hyperbolic flows that mix but at an arbitrarily
slow rate~\cite{Ruelle83,Pollicott85}.   We note that so far, exponential 
decay of correlations has been proved only in very special 
cases~\cite{Dolgopyat98a,Liverani04,Pollicott99}.

In previous work~\cite{Mapp}, we extended Dolgopyat's results to a class
of nonuniformly hyperbolic flows.   These flows are the continuous time
analogue of a class of discrete time nonuniformly hyperbolic systems that are
known, by the results of Young~\cite{Young98}, to have exponential decay of 
correlations.  In this context, we proved that again the flows typically 
mix faster than any polynomial rate.

In this paper, we consider nonuniformly hyperbolic flows for which the 
analogous class of discrete time system is known, by Young~\cite{Young99}, 
to have polynomial decay of correlations.   We show that the flows typically
have polynomial decay of correlations too, with the same polynomial rate 
(as upper bound).

The general set up is that $T:M\to M$ is a nonuniformly hyperbolic 
diffeomorphism in the 
sense of Young~\cite{Young98} but with a polynomial return time function $r$
as in~\cite{Young99}.
In particular, $T:M\to M$ is modelled by a tower $f:\Delta\to\Delta$ 
constructed over a ``uniformly hyperbolic'' base $Y\subset M$.
The degree of nonuniformity is measured by the return time function
$r:Y\to\Z^+$ to the base.    It is assumed that $\Lambda$ intersects its unstable manifolds in positive Lebesgue measure sets and that
$\int r\,d\mu^u<\infty$ where $\mu^u$ denotes Lebesgue measure on unstable manifolds.   Then there exists a physical (SRB) $T$-invariant ergodic probability
measure $\nu$.

Given a H\"older continuous roof function $h:M\to\R^+$, form the suspension 
$M^h=\{(x,u)\in M\times\R:0\le u\le h(x)\}/\sim$ where $(x,h(x))\sim(Tx,0)$.
The suspension flow $\phi_t:M^h\to M^h$ is given by
$\phi_t(x,u)=(x,u+t)$ computed modulo identifications with
$\phi_t$-invariant ergodic 
probability measure $\nu^h=\nu\times{\rm lebesgue}/\int_Mh\,d\nu$.

Suppose that the return time function satisfies the 
polynomial tails condition 
\[
\mu^u(y\in Y:r(y)>n)=O(1/n^{\beta+1}), \quad \beta>0.
\]
Under this assumption, 
Young~\cite{Young99} obtained the decay rate
\begin{align} \label{eq-maps}
\SMALL \int_M v\,w\circ T^n\,d\nu -\int_M v\,d\nu \int_M w\,d\nu=O(1/n^{\beta}),
\end{align}
for the discrete time dynamics and H\"older observables
$v,w:M\to\R$.   We prove that typically the 
decay of correlations for the suspension flow satisfies
\begin{align} \label{eq-flows}
\SMALL\rho_{v,w}(t)=\int_{M^h} v\,w\circ\phi_t\,d\nu^h -
\int_{M^h} v\,d\nu^h \int_{M^h} w\,d\nu^h =O(1/t^\beta),
\end{align}
provided $v,w:M^h\to\R$ are sufficiently regular.

\begin{rmk}    \label{rmk-sharp}
(i) 
In certain situations, estimate~\eqref{eq-maps} is 
sharp~\cite{Gouezel04a,Hu04,Sarig02} and it seems likely that
estimate~\eqref{eq-flows}
 is also sharp in the generality that it is proved.
\\[.75ex]
\noindent(ii)   As in~\cite{Dolgopyat98b,Mapp},
the results in this paper hold only for observables that are sufficiently
smooth in the flow direction.
In particular, our results do not apply
to the position variable in the Lorentz flow examples below.
\\[.75ex]
\noindent(iii) The approach in this paper works for general decay rates of 
$\mu^u(y\in Y:r(y)>n)$, see Sections~\ref{sec-trunc} and~\ref{sec-exp}.    
We note that the calculations are considerably simpler in the special case
$\mu^u(y\in Y:r(y)>n)=O(1/n^{\beta+1})$, $\beta>1$.   
\\[.75ex]
\noindent(iv) 
Throughout this paper, we require that the roof function $h$ is bounded below
away from zero.   
(Such an assumption was not required in~\cite{Mapp}.)
Intuitively, one expects that the violation of this
condition may actually accelerate mixing (for example in the case of
Sina{\u\i} billiards with cusps, see~\cite[p.~15]{ChernovDolgopyatapp}
or~\cite[Section~5.6]{ChernovYoung00}).   
\end{rmk}

Our methods apply also to suspension flows with unbounded
roof functions.  We assume an exponential tails condition
$\mu^u(y\in Y:r(y)>n)=O(\gamma^n)$, $\gamma\in(0,1)$ for the return time $r$,
and consider an unboundedness assumption of the type
\[
\mu^u(x\in M:h(x)> n)=O(1/n^{\beta+1}),\quad \beta>0,
\]
for the roof function $h$.
(For technical reasons, the actual assumption is slightly more complicated,
see Section~\ref{sec-main-c} for the precise statement.)
We prove that typically, for $\epsilon>0$ arbitrarily small 
\begin{align} \label{eq-eps}
\rho_{v,w}(t)=O(1/t^{\beta-\epsilon}).
\end{align}

\begin{rmk}    \label{rmk-sharp2}
In future work, we expect to improve estimate~\eqref{eq-eps} to 
\begin{align} \label{eq-unbounded}
\rho_{v,w}(t)=O((\ln t)^{\beta+1}/t^\beta).
\end{align}
Indeed we obtain this improved estimate for nonuniformly expanding 
semiflows, see Theorem~\ref{thm-unbounded}.

It is possible that the logarithmic factor in~\eqref{eq-unbounded} is an 
artifact of
our method, but it seems more likely that additional assumptions on $r$ and $h$
jointly are required to remove it.   
This is discussed in Section~\ref{sec-log}.   Such conditions
may be harder to verify in practice.

In contrast, the logarithmic factors in Examples~\ref{ex-semi} 
and~\ref{ex-curv} below are
inherited from the tail estimates for $r$ 
in~\cite{ChernovZhang05,ChernovZhang05b,Markarian04}.
\end{rmk}

We now list some applications of our results.   A good source of examples
are provided by billiards and the associated Lorentz flows~\cite{ChernovDolgopyatapp,ChernovYoung00}.  However, it should be emphasized that the results
apply generally to nonuniformly hyperbolic systems modelled by Young towers.
In particular, whereas it is likely that decay of correlations hold at the
specified rate for {\em all} Lorentz flows in the examples below, it is 
well-known even in the uniformly hyperbolic context, and hence certainly in the 
generality of this paper,
that positive results can be expected only for {\em typical} flows.
As in~\cite{Mapp}, it suffices that any four periodic orbits intersecting
the base $Y$ have periods satisfying a Diophantine-type condition,
see Corollary~\ref{cor-data}.
(In the uniformly hyperbolic case, it suffices to consider any pair of periodic 
orbits~\cite{Dolgopyat98b}.)

\begin{example}[Intermittency-type semiflows]
Various authors including~\cite{Hu04,LiveraniSaussolVaienti99,Young99}
have studied 
intermittency (Pomeau-Manneville) maps of the type $T:[0,1]\to[0,1]$
given by
\[
Tx=\left\{ \begin{array}{cc}
x(1+2^{\alpha}x^{\alpha})& 0\le x<\frac{1}{2}\\
2x-1 &\frac{1}{2}\le x<1
\end{array}\right.
\]
for $0<\alpha<1$,  where there is an indifferent fixed point at $0$.
There is a unique absolutely continuous ergodic invariant
probability measure $\nu$ and for $\eta>0$ there is a constant $C$ such that
\[
\SMALL \bigl|\int_{[0,1]} v\, w\circ T^n\, d\nu-\int_{[0,1]} v\,d\nu\int_{[0,1]} w\,d\nu\bigr|
\le C\|v\|_{C^\eta} |w|_\infty /n^\beta, \qquad \beta=\frac{1}{\alpha}-1,
\]
for all $v\in C^\eta([0,1])$, $w\in L^\infty([0,1])$, $n\ge1$.
Hu~\cite{Hu04} showed that this rate is optimal.
Further, the upper bound $O(1/n^{\beta})$ was obtained in~\cite{Young99}
via the the construction of a Young tower with tail decay rate
$1/n^{\beta+1}$.

Now construct the suspension semiflow $\phi_t:[0,1]^h\to [0,1]^h$ where 
$h:[0,1]\to\R^+$ is a H\"older continuous roof function.
Assume that $v:[0,1]^h\to\R$ is $C^\infty$ along the flow direction
with $C^\eta$ derivatives for some $\eta>0$.   
For typical roof functions $h$, it follows from Theorem~\ref{thm-exp} below 
that there exists a constant $C_v$ such that
\[
|\rho_{v,w}(t)|\le C_v |w|_\infty /t^\beta,
\]
for all $w\in L^\infty([0,1]^h)$, $t>0$.
\end{example}

\begin{example}[Semi-dispersing Lorentz flows]
\label{ex-semi}
Chernov \& Zhang~\cite{ChernovZhang05} consider 
a class of semidispersive billiards with tables of the form
$R-\{B_1\cup \dots\cup B_r\}$ where $R$ is a rectangle and 
$B_1,\dots,B_r\subset\operatorname{Int}R$ are disjoint strictly convex 
obstacles with $C^3$ boundaries
(see~\cite[Figure~2(a)]{ChernovZhang05}).
Building upon ideas of Markarian~\cite{Markarian04}, it is shown
in~\cite{ChernovZhang05} that
the correlation function for the billiard map (for H\"older observables)
decays as $O((\ln n)^2/n)$.
A byproduct of the proof (see~\cite[Section~3]{ChernovZhang05}) is the 
existence of a Young tower with tails decaying as $O((\ln n)^2/n^2)$.     
Hence, it follows from Theorem~\ref{thm-hyp} below
(noting Remark~\ref{rmk-h})
that typically the corresponding Lorenz flows have decay rates 
$\rho_{v,w}(t)=O((\ln t)^2/t)$ for observables $v,w$ sufficiently
smooth in the flow direction.
\end{example}

\begin{example}[Dispersing Lorentz flows with vanishing curvature]
\label{ex-curv}

Chernov \& Zhang~\cite{ChernovZhang05b} study a
class of dispersing billiards where the billiard table
has smooth strictly convex boundary with nonvanishing curvature,
except that the curvature vanishes at two points.
Moreover, it is assumed that there is a periodic orbit that runs between the 
two flat points, and that the boundary near these flat points has the form
$\pm(1+|x|^b)$ for some $b>2$.
The correlation function for the billiard map 
decays as $O((\ln n)^{\beta+1}/n^\beta)$ 
where $\beta=(b+2)/(b-2)\in(1,\infty)$.
Again, a byproduct of the proof is the 
existence of a Young tower with tails decaying as $O((\ln n)^{\beta+1}/n^{\beta+1})$.     
Hence, it follows from Theorem~\ref{thm-hyp}
that typically the corresponding Lorenz flows have decay rates $O((\ln t)^{\beta+1}/t^\beta)$.
\end{example}

\begin{example}[Infinite horizon planar periodic Lorentz gas]
\label{ex-horizon}
The planar periodic Lorenz gas is a class of examples introduced 
by Sina{\u\i}~\cite{Sinai70}.
The billiard map $T:M\to M$ has exponential decay of correlations, as shown
by Young~\cite{Young98} in the finite horizon case and Chernov~\cite{Chernov99}
in the infinite horizon case.  In both cases, the map is modelled by a
Young tower with exponential tails.

In the finite horizon case, Chernov~\cite{Chernovsub} has recently proved that
correlations for the Lorentz flow decay at least stretched exponentially.
(Previously~\cite{Mapp} showed that the decay is
typically faster than any polynomial rate.)
For the infinite horizon case, it is widely expected that the decay
rate is $1/t$, see~\cite{FriedmanMartin88}.
A calculation shows that $\mu(x\in M:h(x)>n)=O(1/n^2)$, and
it follows from Theorem~\ref{thm-unbounded} below (noting Remark~\ref{rmk-eps})
that typically the Lorentz flow has decay rate $O(t^{1-\epsilon})$ where
$\epsilon>0$ is arbitrarily small.
\end{example}

The remainder of the paper is organised as follows.
In Section~\ref{sec-main}, we state the main results first in the simpler
context of nonuniformly expanding semiflows, and then for 
nonuniformly hyperbolic flows.
Also, we present an outline of the strategy of the proof.
The proof for nonuniformly expanding semiflows is carried out
in Section~\ref{sec-trunc} and~\ref{sec-exp}.
The modifications for nonuniformly hyperbolic flows are described in
Section~\ref{sec-hyp}.
The modifications for (semi)flows with unbounded roof function
are described in Section~\ref{sec-unbounded}.

\section{Statement of the main results}
\label{sec-main}

In this section, we state our main results.
In Subsection~\ref{sec-main-a}, we consider the technically simpler case of 
nonuniformly expanding semiflows.
In Subsection~\ref{sec-main-b}, we consider nonuniformly hyperbolic flows.  
The case of unbounded roof function is discussed in Subsection~\ref{sec-main-c}.
In Subsection~\ref{sec-main-d}, we describe the strategy of the proof, 
focusing for simplicity on the result in Subsection~\ref{sec-main-a}.

\subsection{Nonuniformly expanding semiflows}
\label{sec-main-a}

Let $(X,d)$ be a locally compact separable bounded metric space with
Borel probability measure $\mu_0$ and let $T:X\to X$ be a nonsingular
transformation for which $\mu_0$ is ergodic.
Let $Y\subset X$ be a measurable subset with $\mu_0(Y)>0$, and
let $\{Y_j\}$ be an at most countable measurable partition
of $Y$ with $\mu_0(Y_j)>0$.    We suppose that there is an $L^1$
{\em return time} function $r:Y\to\Z^+$, constant on each $Y_j$ with
value $r(j)\ge1$, and constants $\lambda>1$, $\eta\in(0,1)$, $C\ge1$
such that for each $j\ge1$,
\begin{itemize}
\item[(1)] $F=T^{r(j)}:Y_j\to Y$ is a bijection.
\item[(2)] $d(Fx,Fy)\ge \lambda d(x,y)$ for all $x,y\in Y_j$.
\item[(3)] $d(T^\ell x,T^\ell y)\le Cd(Fx,Fy)$ for all $x,y\in Y_j$,
$0\le \ell <r(j)$.
\item[(4)] $g_j=\frac{d(\mu_0|{Y_j}\circ F^{-1})}{d\mu_0|_Y}$
satisfies $|\log g_j(x)-\log g_j(y)|\le Cd(x,y)^\eta$ for all
$x,y\in Y$.
\end{itemize}
Such a map $T:X\to X$ is called {\em nonuniformly expanding}.
There is a unique $T$-invariant probability measure $\nu$ equivalent
to $\mu_0$ (see for example~\cite[Theorem~1]{Young99}).

\begin{rmk}   Discarding sets of zero measure, we have assumed without loss
that the induced map $F:Y\to Y$ is defined everywhere on $Y$.
This simplifies the formulation below of certain hypotheses involving
periodic points.
\end{rmk}

Let $h:X\to\R^+$ be a roof function such that for all $j\ge1$,
\begin{itemize}
\item[(5)] $h,\frac1h\in L^\infty(X)$ and
$|h(x)-h(y)|\le Cd(x,y)^\eta$ for all $x,y\in T^\ell Y_j$, $0\le\ell< r(j)$.
\end{itemize}
Define the suspension semiflow $T_t:X^h\to X^h$ with invariant ergodic 
measure $\nu^h$ as in the introduction.
Let $\rho_{v,w}(t)$ denote the correlation function 
corresponding to observables $v,w:X^h\to\R$.

For $m\ge1$, $\eta>0$, let $C^{m,\eta}(X^h)$ consist of those $v:X^h\to\R$ for which
$\|v\|_{m,\eta}=\|v\|_\eta+\|\partial_tv\|_\eta
+\dots+\|\partial_t^mv\|_\eta<\infty$,
where $\partial_t$ denotes the derivative in the flow direction and
\[
\|v\|_\eta = |v|_\infty+\sup_{(x,u)\neq(y,u)}|v(x,u)-v(y,u)|/d(x,y)^\eta.
\]

Suppose that  $Z\subset Y$ is a finite union of partition elements $Y_j$.
Let $p\in Z$ be a periodic point for $F:Y\to Y$ such that $F^ip\in Z$ for
all $i\ge1$.   We associate to $p$ the triple
$(\tau,d,q)\in\R^+\times\Z^+\times\Z^+$
where $\tau$ is the period of $p$ under the semiflow $T_t$, $d$ is the period
under the map $T$, and $q$ is the period under the induced map $F$ (so $d=\sum_{i=0}^{q-1}r(F^ip)$ and
$\tau=\sum_{i=0}^{d-1}h(T^ip)$).
Let $\mathcal{T}_Z$ denote the set of such triples.

\begin{thm}  \label{thm-exp}
Let $T:X\to X$ be a nonuniformly expanding map and $h:X\to\R^+$ a roof function
satisfying properties (1)--(5).
Assume that $\mu_0(y\in Y:r(y)>n)=O((\ln n)^\gamma n^{-(\beta+1)})$, 
for some $\beta>0$, $\gamma\ge0$.
Let $Z\subset Y$ be a finite union of partition elements $Y_j$.

Suppose that there do {\bf not} exist constants $C,m\ge1$ such that
\[
|\rho_{v,w}(t)|\le C\|v\|_{m,\eta}|w|_\infty (\ln t)^{\gamma}t^{-\beta},
\]
for all $t>0$, $v\in C^{m,\eta}(X^h)$, $w\in L^\infty(X^h)$.

Then there exist sequences $b_k\in\R$ with $|b_k|\to\infty$,
and $\omega_k,\,\varphi_k\in[0,2\pi)$;
and constants $\alpha>0$ arbitrarily large, $C,\beta_0\ge1$; such that
\begin{align} \label{eq-data}
\dist(b_kn_k\tau+\omega_kn_kd+q\varphi_k,2\pi\Z)\le Cq |b_k|^{-\alpha},
\end{align}
 for all $k\ge1$ and all $(\tau,d,q)\in\mathcal{T}_Z$, where $n_k=[\beta_0\ln|b_k|]$.
\end{thm}

\begin{rmk} \label{rmk-h}
It is easy to relax the condition that $1/h\in L^\infty$
to the requirement that there exists an $n_0\ge0$ such 
that $1/h_{n_0}\in L^\infty$ where $h_{n_0}=h+h\circ T+\dots+h\circ T^{n_0-1}$.
In Example~\ref{ex-semi}, this condition
is satisfied for $n_0=2$ even though $h$ is arbitrarily small near
the four corner points.
\end{rmk}

\begin{cor} \label{cor-data}
Let $T:X\to X$ be a nonuniformly expanding map and $h:X\to\R^+$ a roof function
satisfying properties (1)--(5).
Assume that $\mu_0(y\in Y:r(y)>n)=O((\ln n)^\gamma n^{-(\beta+1)})$ 
for some $\beta>0$, $\gamma\ge0$.

There exists an integer $m$ with the following property:
Fix four periodic solutions for $T_t:X^h\to X^h$ that each intersect $Y$,
and let $\tau_1,\ldots,\tau_4$ be the periods.
For Lebesgue almost all $(\tau_1,\cdots,\tau_4)\in(\R^+)^4$, 
there exists a constant $C\ge1$ such that
\[
|\rho_{v,w}(t)|\le C\|v\|_{m,\eta}|w|_\infty\,(\ln t)^{\gamma}t^{-\beta},
\]
for all $t>0$, $v\in C^{m,\eta}(X^h)$, $w\in L^\infty(X^h)$.
\end{cor}

\begin{proof}  See~\cite[Corollary~2.4]{Mapp}.
\end{proof}

\begin{rmk} \label{rmk-good}
Similarly, it suffices that there is a sequence of periodic orbits
in $Z$ with {\em good asymptotics} in the sense of~\cite{FMTapp}.
By~\cite{FMTapp}, good asymptotics is an open-dense
condition for smooth systems.   Hence results on
{\em stable} rates of mixing reduce to stability of the partition $\{Y_j\}$.
We do not explore this issue further in this paper.
\end{rmk}

\subsection{Nonuniformly hyperbolic flows}
\label{sec-main-b}

Let $(M,d)$ be a Riemannian manifold. Young~\cite{Young98}
introduced a class of nonuniformly hyperbolic diffeomorphisms $T:M\to M$
(possibly with singularities) with the property that there is an
ergodic $T$-invariant SRB measure $\nu$ for which exponential decay of
correlations holds for H\"older observables.
We refer to~\cite{Young98} for precise definitions, but some of the
notions and notation are required to state our main results.
(The further structure from~\cite{Young98} required for our proofs
is made explicit in Section~\ref{sec-hyp}.)
 In particular, there is a ``uniformly hyperbolic'' subset $Y\subset M$ 
with partition $\{Y_j\}$ and a return time function $r:Y\to \Z^+$ constant 
on partition elements such that, modulo uniformly contracting directions, 
the induced map $F=T^r:Y\to Y$ is nonuniformly expanding.

The statement of our main result is completely analogous to that
of Theorem~\ref{thm-exp}.   Given a roof function
$h:M\to\R^+$, the suspension flow $T_t:M^h\to M^h$ and ergodic measure
$\nu^h$ is defined as before.
Suppose that  $Z\subset Y$ is a finite union of partition elements
$Y_j$.  Again, we define
the set $\mathcal{T}_Z$ consisting of triples $(\tau,d,q)$ corresponding
to periodic orbits for $F:Y\to Y$ lying entirely in $Z$.

\begin{thm}  \label{thm-hyp}
Let $T:M\to M$ be nonuniformly hyperbolic in the sense of Young~\cite{Young98}
with $\mu^u(y\in Y:r(y)>n)=O((\ln n)^\gamma n^{-(\beta+1)})$ 
for some $\beta>0$, $\gamma\ge0$.
Let $h:M\to\R^+$ be a roof function with $h,\frac1h\in L^\infty(M)$ and
$|h(x)-h(y)|\le Cd(x,y)^\eta$ for all $x,y\in T^\ell Y_j$,
$0\le\ell< r(j)$.
Let $Z\subset Y$ be a finite union of partition elements $Y_j$.

Suppose that there do {\bf not} exist constants $C,m\ge1$ such that
\[
|\rho_{v,w}(t)|\le C\|v\|_{m,\eta}\|w\|_{m,\eta}\,(\ln t)^{\gamma}t^{-\beta},
\]
for all $t>0$, $v,w\in C^{m,\eta}(M^h)$.

Then condition~\eqref{eq-data} holds as in Theorem~\ref{thm-exp}.
\end{thm}

\subsection{Unbounded roof functions}
\label{sec-main-c}

Suppose now that $T$ is nonuniformly expanding as in 
Subsection~\ref{sec-main-a}, except that 
the roof function $h:X\to\R^+$ may be unbounded.   Condition (5) is relaxed to
\begin{itemize}
\item[(6)] 
$\frac1h\in L^\infty(X)$ and $|h(x)-h(y)|\le Cd(x,y)^\eta$ for all 
$x,y\in T^\ell Y_j$, $0\le\ell< r(j)$.
\end{itemize}
Also, we require the following technical assumption on the growth of $h$.
Let $X(n)=\bigcup\{T^\ell Y_j:\|h1_{T^\ell Y_j}\|_{C^\eta}\ge n\}$.
\begin{itemize}
\item[(7)]  $\mu_0(X(n))\le Cn^{-(\beta+1)}$.
\end{itemize}

\begin{thm} \label{thm-unbounded}
Let $T:X\to X$ be a nonuniformly expanding map satisfying
properties (1)--(4) and $h:X\to\R^+$ a roof function
satisfying properties (6), (7) where $\beta>0$.
Assume that $\mu_0(y\in Y:r(y)>n)=O(e^{-cn})$ for some $c>0$.
Let $Z\subset Y$ be a finite union of partition elements $Y_j$.

Suppose that there do {\bf not} exist constants $C,m\ge1$ such that
\[
|\rho_{v,w}(t)|\le C\|v\|_{m,\eta}|w|_\infty (\ln t)^{\beta+1}t^{-\beta},
\]
for all $t>0$, $v\in C^{m,\eta}(X^h)$, $w\in L^\infty(X^h)$.

Then condition~\eqref{eq-data} holds as in Theorem~\ref{thm-exp}.
\end{thm}

\begin{rmk}
The methods in this paper can handle general decay rates for
$\mu_0(r>n)$ and $\mu_0(X(n))$.   In the absence of motivating
examples, we do not consider this generality.   Again,
a joint estimate of these quantities,
if available in a specific application, might lead to improved results,
as in Remark~\ref{rmk-sharp2}.
\end{rmk}

\begin{rmk} \label{rmk-eps}
An analogous result holds for $T:M\to M$ nonuniformly hyperbolic,
except that presently we obtain the weaker decay rate~\eqref{eq-eps}.
\end{rmk}

\subsection{Strategy of the proof}
\label{sec-main-d}

There are a number of steps in proving Theorem~\ref{thm-exp}.

\vspace{-1ex}
\paragraph{Step 1}  We model the nonuniformly expanding map $T:X\to X$ 
by a tower map $f:\Delta\to\Delta$.
Recall that $F=T^r:Y\to Y$ is the induced map.
Define $\Delta=\{(y,\ell)\in Y\times\N: 0\le\ell\le r(y)\}/\sim$ where 
$(y,r(y))\sim(Fy,0)$.
Define the tower map $f:\Delta\to \Delta$ by setting 
$f(y,\ell)=(y,\ell+1)$ computed modulo identifications.
The projection $\pi:\Delta\to X$, $\pi(y,\ell)=T^\ell y$ defines a 
semiconjugacy, $\pi\circ f=T\circ \pi$.

There is a unique invariant ergodic probability measure $\mu_Y$
equivalent to $\mu_0|Y$ for the induced map $F:Y\to Y$.    
Moreover, the density is bounded above (and below) so that $\mu_Y$
inherits the property $\mu_Y(r>n)=O((\ln n)^\gamma n^{-(\beta+1)})$.

We obtain an invariant probability measure on $\Delta$ given by 
$\mu_\Delta=\mu_Y\times\mu_C/\int_Y r\,d\mu_Y$ where 
$\mu_C$ denotes counting measure,
and $\pi:\Delta\to X$ is measure-preserving, carrying $\mu_\Delta$ to $\nu$.
Given $h:X\to\R^+$ H\"older, let $\tilde h=h\circ\pi:\Delta\to\R^+$.
We obtain invariant measures $\nu^h$ and 
$\mu_{\Delta^{\tilde h}}=(\mu_\Delta)^{\tilde h}$ for the suspension flows
on $X^h$ and $\Delta^{\tilde h}$.   The projection $\pi:\Delta\to X$
induces a projection $\pi:\Delta^{\tilde h}\to X^h$ which carries
$\mu_{\Delta^{\tilde h}}$ to $\nu^h$.

If $x,y\in Y$, let $s(x,y)$ be the least integer $n\ge0$ such that
$F^nx,F^ny$ lie in distinct partition elements in $Y$.
If $x,y\in Y_{j}\times\{\ell\}$, then there exist unique $x',y'\in Y_j$
such that $x=f^\ell x'$ and $y=f^\ell y'$.   Set $s(x,y)=s(x',y')$.
For all other pairs $x,y$, set $s(x,y)=0$.
This defines a separation time $s:\Delta\times\Delta\to\N$
and hence a metric $d_\theta(x,y)=\theta^{s(x,y)}$ on $\Delta$.
Let $F_\theta(\Delta)$ denote the Banach space of Lipschitz functions 
$v:\Delta\to\R$ with norm $\|v\|_\theta=|v|_\infty+|v|_\theta$
where $|v|_\theta=\sup_{x\neq y}|v(x)-v(y)|/d_\theta(x,y)$.
We can choose $\theta\in(0,1)$ so that
$v\circ\pi\in F_\theta(\Delta)$ for all $v\in C^\eta(X)$.
It follows that $v\circ\pi\in F_{m,\theta}(\Delta^{\tilde h})$ 
for all $v\in C^{m,\eta}(X^h)$ where $\tilde h=h\circ\pi$ and
$F_{m,\theta}(\Delta^{\tilde h})$ is defined in the obvious way.

Hence, we may reduce to the situation where $f_t:\Delta^h\to\Delta^h$
is a suspension flow over a tower map $f:\Delta\to\Delta$
and $h\in F_\theta(\Delta)$.   It suffices to consider decay of correlations
for observables $v\in F_{m,\theta}(\Delta^h)$, $w\in L^\infty(\Delta^h)$.

\vspace{-1ex}
\paragraph{Step 2}  We truncate the return time function $r$
so that $r\le N$.   This produces an error 
$O\bigl((\ln N)^\gamma N^{-\beta}+t(\ln N)^\gamma N^{-(\beta+1)}\bigr)$ and 
reduces the problem to a suspension 
flow over a truncated tower $f':\Delta'\to\Delta'$ with bounded return 
time $r'$.
(All constructions from $F:Y\to Y$, $\mu_Y$, $r$ and $h$, are repeated with
$r$ replaced by $r'$.)

\vspace{-1ex}
\paragraph{Step 3}   
Let $\rho'(t)$ denote the correlation function on $(\Delta')^h$ and
let $\hat\rho(s)$ denote the Laplace transform of $\rho'(t)$.
Modulo an analytic term,
\[
\hat\rho(s)\sim\sum_{n\ge1}\int_{\Delta'}
e^{-sh_n\circ\pi}v_s\,w_s\circ (f')^n
\, d\mu_{\Delta'}, 
\]
where
$v_s(x)=\int_0^{h(x)}e^{su}\,v(x,u)\,du$ and
$w_s(x)=\int_0^{h(x)}e^{-su}\,w(x,u)\,du$.
Decay of $\rho'(t)$ reduces to analyticity properties of $\hat\rho(s)$.

\vspace{-1ex}
\paragraph{Step 4}
Let $L:L^1(\Delta')\to L^1(\Delta')$ be the 
(Perron-Frobenius) transfer operator for
$f':\Delta'\to\Delta'$
(so $\int_{\Delta'}v\,w\circ f'\,d\mu_{\Delta'}
=\int_{\Delta'}Lv\,w\,d\mu_{\Delta'}$
for $v\in L^1(\Delta')$, $w\in L^\infty(\Delta')$).
For $s\in\C$, define
the twisted transfer operator $L_s$ by $L_sv=L(e^{sh}v)$.  Then 
\[
\hat\rho(s)\sim\sum_{n\ge1}\int_{\Delta'}L_{-s}^nv_s\,w_s
\, d\mu_{\Delta'}.
\]
Via the technique of operator renewal sequences, 
estimates on $|L_s^nv|_1$ for $v\in F_{m,\theta}((\Delta')^h)$
are related to estimates for the transfer operator of the (fixed) induced map
$F:Y\to Y$.   

\vspace{-1ex}
\paragraph{Step 5} Choosing $N=N(t)$ appropriately, the estimates in Steps~2 and~4
yield the required result.

\vspace{1ex}
Steps~1 and~3 are standard.
See for example~\cite[Section~4.1]{Mapp}
for Step~1, and~\cite{Dolgopyat98a,Pollicott85} or 
specifically~\cite[Section~10]{Dolgopyat98b} for Step~3.
The truncation in Step~2 is the main new idea in this paper and is 
carried out in Section~\ref{sec-trunc}.
In Section~\ref{sec-exp}, we carry out Step~4
following~\cite{Mapp} but keeping careful track of the dependence of
estimates on $N$.
We then specify $N=N(t)$ to obtain the final result.

\begin{rmk}
(a)
In Step~1, the subset $Y\subset X$ is identified with the base 
$\{\ell=0\}$ of the tower $\Delta$.   The induced map $F:Y\to Y$ becomes a 
first return map for $f:\Delta\to\Delta$.  In particular, $f:\Delta\to\Delta$
is Markov, even though no such assumption is made on $T:X\to X$.
\\[.75ex]
(b) The induced map $F:Y\to Y$ is a full shift on a countable alphabet
with good distortion properties (guaranteed by condition~(4) in
Subsection~\ref{sec-main-a}).
Such maps are often called {\em Gibbs-Markov} and are studied extensively
in~\cite{Aaronson}.
\\[.75ex]
(c)  
Note that $h$ is unchanged in Step~2, except that it is 
restricted to $\Delta'$.   Similarly for $v$, $w$, except that a 
further approximation is required to ensure that $v$ remains inside 
$F_{m,\theta}(\Delta^h)$, see Section~\ref{sec-buffer}.
\\[.75ex]
(d)  
The truncation in Step~2 seems at first sight to make uncontrollable
changes to $f:\Delta\to\Delta$ and hence to the suspension flow on
$\Delta^h$.  However, it should be noted that the induced map $F:Y\to Y$
and the $F$-invariant measure $\mu_Y$ are unchanged by the truncation.  
The techniques in~\cite{Mapp}
based on operator renewal sequences~\cite{Sarig02,Gouezel04a,BHM05}
are hence well-suited to this situation, see Section~\ref{sec-exp}.
\end{rmk}

\section{Truncation of the roof function}
\label{sec-trunc}

Let $f:\Delta\to\Delta$ be a tower map, modelling the
underlying nonuniformly expanding map $T:X\to X$,
as discussed in Section~\ref{sec-main-d}.
Let $h:\Delta\to\R^+$ be a Lipschitz roof function and
let $f_t:\Delta^h\to \Delta^h$ be the suspension flow.
Recall that $\Delta$ is itself a discrete suspension over the induced map
$F:Y\to Y$ with ergodic invariant probability measure $\mu_Y$ and return time
$r:Y\to\Z^+$.
The tower map $f:\Delta\to\Delta$ has an ergodic invariant probability
measure $\mu_\Delta=\mu_Y\times\text{counting}/\bar r$ where 
$\bar r=\int_Y r\,d\mu_Y$.
Similarly, $f_t:\Delta^h\to \Delta^h$ has an ergodic invariant probability
measure $\mu_{\Delta^h}=\mu_\Delta\times\text{lebesgue}/\bar h$
where $\bar h=\int_\Delta h\,d\mu_\Delta$.

For fixed $N\ge1$, we define the truncated return time function 
$r'=\min\{r,N\}:Y\to\Z^+$.
Then we form the truncated tower map $f':\Delta'\to\Delta'$ over $Y$ with 
measure $\mu_{\Delta'}=\mu_Y\times\text{counting}/\bar{r'}$.
Restricting $h$ to $\Delta'$, we obtain the truncated suspension flow
$f_t':(\Delta')^h\to (\Delta')^h$ with measure $\mu_{(\Delta')^h}=
\mu_{\Delta'}\times\text{lebesgue}/\bar h'$
where $\bar h'=\int_{\Delta'} h\,d\mu_{\Delta'}$.

Write $\Delta=\Delta_{\rm left}\dot\cup \Delta_{\rm right}$
where
\[
\Delta_{\rm left}=\{(y,\ell)\in\Delta:r(y)<N\}, \qquad
\Delta_{\rm right}=\{(y,\ell)\in\Delta:r(y)\ge N\}.
\]

\begin{prop} \label{prop-r}
\begin{itemize}
\item[(i)] $\bar r-\bar{r'}=\sum_{n>N}\mu_Y(r\ge n)$.
\item[(ii)] 
$\mu_\Delta(\Delta_{\rm right})=(1/\bar r)\{N\mu_Y(r\ge N)+\sum_{n>N}\mu_Y(r\ge n)\}$.
\end{itemize}
\end{prop}

\begin{proof} 
This is a standard computation.
\end{proof}

\begin{prop} \label{prop-Ek}
For $k\ge1$, define
\[
E_k=\{x\in\Delta:f^jx\in\Delta_{\rm right} \enspace\text{for at least
one $j\in\{0,1,\dots,k$\}}\}.
\]
Then
$\mu_\Delta(E_k)\le (1/\bar r)\{\sum_{n>N}\mu_Y(r\ge n)\,+\,
(N+k)\mu_Y(r\ge N)\}$.
\end{prop}

\begin{proof}
Write $E_k$ as the disjoint union $E_k=
\bigcup_{j=0}^k G_j$ where
\[
G_j=\{
f^ix\in\Delta_{\rm left}\enspace\text{for $i\in\{0,1,\dots,j-1$\}}
\enspace\text{and}\enspace
f^jx\in\Delta_{\rm right} 
\}.
\]
In particular $\mu_\Delta(G_0)=\mu_\Delta(\Delta_{\rm right})$.
For $j\ge1$, it follows from the definition that if $x\in G_j$, then
$f^jx\in\Delta_{\rm right}\cap Y$ (the base of the tower).
Hence $\mu_\Delta(G_j)\le
\mu_\Delta(f^{-j}(\Delta_{\rm right}\cap Y))
=
\mu_\Delta(\Delta_{\rm right}\cap Y)=(1/\bar r)\mu_Y(r\ge N)$.
\end{proof}

For notational convenience, we write $\Omega=\Delta^h$ and
$\Omega'=(\Delta')^h$ throughout the remainder of this section.
Throughout the paper
$C$ denotes a universal constant, varying from line to line, dependent only
on the suspension semiflow $T_t:X^h\to X^h$ and the regularity exponents 
$m$, $\eta$.

\begin{lemma} \label{lem-trunc}
Suppose that $h,\frac1h:\Delta\to\R^+$, $v,w:\Omega\to\R$ all lie in 
$L^\infty$.  Let 
\begin{align*}
\rho(t) & =\int_\Omega v\,w\circ f_t\,d\mu_\Omega -\int_\Omega v\,d\mu_\Omega \int_\Omega w\,d\mu_\Omega,
\\
\rho'(t) & =\int_{\Omega'} v\,w\circ f'_t\,d\mu_{\Omega'} -\int_{\Omega'} v\,d\mu_{\Omega'} \int_{\Omega'} w\,d\mu_{\Omega'}.
\end{align*}
Then there exists $N_0$, $t_0$ (depending only on $\mu$, $r$ and $h$) 
such that for all $N\ge N_0$, $t\ge t_0$,
\[
\SMALL|\rho(t)-\rho'(t)|\le 
C|v|_\infty|w|_\infty \{\sum_{n>N}\mu_Y(r\ge n)\,+\, (N+t)\mu_Y(r\ge N)\}.
\]
\end{lemma}

\begin{proof}
We choose $N\ge N_0$ sufficiently large that
$1/\bar r'\le 2/\bar r$, $1/\bar h'\le 2/\bar h$. 
It follows that
$ \frac{1}{\bar r'}-\frac{1}{\bar r}\le \frac{2}{\bar r^2}(\bar r-\bar r')$
and $
|\frac{1}{\bar h'}-\frac{1}{\bar h}|\le \frac{2}{\bar h^2}|\bar h-\bar h'|$.
Further, 
$|\bar h-\bar h'|\le 4|h|_\infty(\bar r-\bar r')/\bar r$.
By Proposition~\ref{prop-r}(i),
\begin{align} \label{eq-frac}
 \frac{1}{\bar r'}-\frac{1}{\bar r}\le C \sum_{n>N}\mu_Y(r\ge n), \quad
\Bigl|\frac{1}{\bar h'}-\frac{1}{\bar h}\Bigr|\le  \sum_{n>N}\mu_Y(r\ge n).
\end{align}

Let $A=\int_\Omega v\,w\circ f_t\,d\mu_\Omega$, 
$A'=\int_{\Omega'} v\,w\circ f'_t\,d\mu_{\Omega'}$.
By definition, $A=(1/\bar h)(1/\bar r)B$, 
$A'=(1/\bar h')(1/\bar{r'})B'$ where
\begin{align*}
B &= \int_Y\sum_{\ell=0}^{r(y)-1}\int_0^{h(y,\ell)}
v(y,\ell,u)\,w\circ f_t(y,\ell,u)\,du\,d\mu_Y, \\ 
B' & =\int_Y\sum_{\ell=0}^{r'(y)-1}\int_0^{h(y,\ell)}
v(y,\ell,u)\,w\circ f'_t(y,\ell,u)\,du\,d\mu_Y.
\end{align*}
Note that $B=\bar r\int_\Delta\int_0^h v\,w\circ f_t\,du\,d\mu_\Delta$
so that $|B|\le \bar r\bar h|v|_\infty|w|_\infty$.    
Hence by~\eqref{eq-frac},
\[
\SMALL |A-A'|\le C\{|v|_\infty|w|_\infty \sum_{n>N}\mu_Y(r\ge n)+|B-B'|\}.
\]

Write $|B-B'|\le I+II$
where 
\[
I=|v|_\infty|w|_\infty\int_Y\sum_{\ell=N}^{r(y)-1}h(y,\ell)1_{\{r(y)>N\}}d\mu_Y
\le C|v|_\infty|w|_\infty \sum_{n>N}\mu_Y(r\ge n),
\]
and
\begin{align*}
II & =
|v|_\infty\int_Y\sum_{\ell=0}^{r'(y)-1}\int_0^{h(y,\ell)}|w\circ f_t(y,\ell,u)-
w\circ f'_t(y,\ell,u)|du\,d\mu_Y 
\\ &
= |v|_\infty\,\bar{r'}\int_{\Delta'}\int_0^{h(x)}|w\circ f_t(x,u)-
w\circ f'_t(x,u)|du\,d\mu_{\Delta'} 
\\ &
= |v|_\infty\,\bar{r}\int_{\Delta'}\int_0^{h(x)}|w\circ f_t(x,u)-
w\circ f'_t(x,u)|du\,d\mu_{\Delta}.
\end{align*}
(Starting from the last expression, we are regarding $\Delta'$ as 
a subset of $\Delta$; for measurable sets
$E\subset\Delta'\subset\Delta$ note that 
$\bar{r'}\mu_{\Delta'}(E)=\bar{r}\mu_{\Delta}(E)$.)

Now $f_t(x,u)=f_t'(x,u)$ provided $N$ is sufficiently large
that $f_s(x,u)$ lies in the part of the suspension 
over $\Delta_{\rm left}$ for $s\in[0,t]$.
Note also that the flow reaches the roof at most
$t|\frac1h|_\infty+1$ times by time $t$ so
it suffices that $f^jx\in\Delta_{\rm left}$ for
$0\le j\le {[t|\frac1h|_\infty]+2}$.
Hence
\[
II\le C|v|_\infty |w|_\infty \mu_\Delta(E_k),
\]
where $k= [t|\frac1h|_\infty]+2$.
By Proposition~\ref{prop-Ek},
\[
\SMALL
II\le C|v|_\infty |w|_\infty\{\sum_{n>N}\mu_Y(r\ge n)\,+\,
(N+k)\mu_Y(r\ge N)\},
\]
and so 
\[
\SMALL
|A-A'|\le C|v|_\infty|w|_\infty
\{\sum_{n>N}\mu_Y(r\ge n)\, + \, (N+t) \mu_Y(r\ge N)\}.
\]
A similar (but simpler) calculation shows that 
\[
\SMALL \bigl|\int_\Omega v\,d\mu_\Omega \int_\Omega w\,d\mu_\Omega-
\int_{\Omega'} v\,d\mu_{\Omega'} \int_{\Omega'} w\,d\mu_{\Omega'}\bigr|
\le C|v|_\infty|w|_\infty \sum_{n>N}\mu_Y(r\ge n),
\]
and the result follows.
\end{proof}

\begin{rmk}  \label{rmk-billiards} 
If $\mu_Y(r\ge n)=O((\ln n)^{\gamma}n^{-(\beta+1)})$, $\beta>0$, $\gamma\ge0$ 
then
\[
\sum_{n>N}\mu_Y(r\ge n)\le C 
\sum_{n\ge N}\frac{(\ln n)^\gamma}{n^{\frac12\beta}}\frac{1}{n^{\frac12\beta+1}}
\le C 
\frac{(\ln N)^{\gamma}}{N^{\frac12\beta}}
\sum_{n\ge N} \frac{1}{n^{\frac12\beta+1}}
\le C \frac{(\ln N)^{\gamma}}{N^\beta}.
\]
Hence
\[
|\rho(t)-\rho'(t)|\le C|v|_\infty|w|_\infty 
\{(\ln N)^{\gamma}N^{-\beta}+ t (\ln N)^{\gamma}N^{-(\beta+1)}\}.
\]
\end{rmk}

\subsection{Regularity of observables}
\label{sec-buffer}

The original observable
$v:\Delta^h\to\R$ is assumed to be smooth in the flow direction,
but after restriction to $(\Delta')^h$ this condition is typically violated
since the identifications are different.   Namely, we now have
\[
(y,r'(y),h(y))\sim (Fy,0,0),
\]
whereas $v$ is smooth respect to the old identifications 
\[
(y,r(y),h(y))\sim (Fy,0,0).
\]
This problem is resolved by using the top level of the truncated tower 
as a {\em buffer}.   That is, we modify $v$ on the strip
$\{(y,N,u):r'(y)=N, u\in[0,h(y,\ell)]\}$ to obtain a new
observable $\tilde v$ that is as regular in the flow direction
on $(\Delta')^h$ as $v$ was on $\Delta^h$.   
Since $v$ is bounded and $h$ is bounded below,
we can make this modification in such a way that 
$\|\tilde v\|_{m,\eta}\le C\|v\|_{m,\eta}$.
The resulting error in the correlation function
is at most $C|v|_\infty|w|_\infty\mu_Y(r=N)$ which is smaller than 
the error in Lemma~\ref{lem-trunc}.
Hence without loss we may suppose that 
the observable $v$ retains its smoothness
in the flow direction when restricted to $(\Delta')^h$.

\section{Decay for nonuniformly expanding semiflows}
\label{sec-exp}

In this section, we complete the proof of Theorem~\ref{thm-exp}.

Define the induced roof function $H:Y\to\R$ by
$H(y)=\sum_{\ell=0}^{r(y)-1}h\circ f^j(y)$.
For $b\in\R, \omega\in[0,2\pi)$, we define 
$M_{b,\omega}:L^\infty(Y)\to L^\infty(Y)$,
\[
M_{b,\omega}v=e^{-ibH}e^{-i\omega r}v\circ F.
\]

\begin{defn}  A subset $Z_0\subset Y$ is a {\em finite subsystem} of $Y$
if $Z_0=\cap_{n\ge1} F^{-n}Z$ where $Z$ is the union of finitely many
elements from the partition $\{Y_j\}$.
(Note that $F|_{Z_0}:Z_0\to Z_0$ is a a full one-sided
shift on finitely many symbols.)
\end{defn}

\begin{defn} \label{def-M}
We say that $M_{b,\omega}$ has an {\em approximate eigenfunction} on a subset
$Z\subset Y$
if there exist constants $\alpha>0$ arbitrarily large, $\beta_0>0$ and $C\ge1$,
and sequences $|b_k|\to\infty$, $\omega_k\in[0,2\pi)$,
$\varphi_k\in [0,2\pi)$, $u_k\in F_\theta(Y)$ with $|u_k|\equiv1$, such that
setting $n_k=[\beta_0\ln |b_k|]$,
\[
|(M_{b_k,\omega_k}^{n_k}u_k)(y)-e^{i\varphi_k}u_k(y)|\le C|b_k|^{-\alpha},
\]
for all $y\in Z$ and all $k\ge1$.
\end{defn}

Define $d_N=\sum_{k=1}^N k \mu_Y(r\ge k)$.
The main result of this section is:
\begin{thm} \label{thm-eigen}
Let $Z_0\subset Y$ be a finite subsystem and suppose that $M_{b,\omega}$ has
no approximate eigenfunctions on $Z_0$.  Choose $N$ sufficiently
large that $r|Z_0\le N$.

Let $d$, $p>0$.  There exists $C$, $m\ge1$, $\epsilon>0$ such that
\[
|\rho'_{v,w}(t)|\le C\|v\|_{m,\theta}|w|_\infty \{d_NN^{1+d}e^{-\epsilon N^{-1}
\ln N\,t}
+(d_NN^d)^{\,p+2}t^{-p}\},
\]
for all $t>0$, $v\in F_{m,\theta}((\Delta')^h)$, $w\in L^\infty((\Delta')^h)$.

\end{thm}

\begin{cor} \label{cor-eigen}
Let $Z_0\subset Y$ be a finite subsystem and suppose that $M_{b,\omega}$ has
no approximate eigenfunctions on $Z_0$.  Choose $N$ sufficiently
large that $r|Z_0\le N$.

If $\mu_Y(r\ge n)=O((\ln n)^\gamma n^{-(\beta+1)})$, $\beta>0$, $\gamma\ge0$, 
then there exist constants $C,m\ge1$ such that
\[
|\rho_{v,w}(t)|\le C\|v\|_{m,\theta}|w|_\infty (\ln t)^{\gamma}
t^{-\beta},
\]
for all $t>0$, $v\in F_{m,\theta}(\Delta^h)$, $w\in L^\infty(\Delta^h)$.
\end{cor}

\begin{proof}
Combining Lemma~\ref{lem-trunc} (specifically Remark~\ref{rmk-billiards})
and Theorem~\ref{thm-eigen}, we obtain
\[
\rho(t)=O\{(\ln N)^\gamma N^{-\beta}+ t(\ln N)^\gamma N^{-(\beta+1)}
+d_N N^{1+d}e^{-\epsilon N^{-1}\ln N\,t} +(d_NN^d)^{\,p+2}t^{-p}\}.
\]
We compute that $d_N\le C$ for $\beta>1$, 
$d_N\le C (\ln N)^{\gamma+1}$ for $\beta=1$, and
$d_N\le C (\ln N)^\gamma N^{1-\beta}$ for $\beta\in(0,1)$.

Set $N=[t/q]$.  Then the first two terms in $\rho(t)$
are $O((\ln t)^\gamma t^{-\beta})$.  The third term
is $O(d_N t^{1+d-\epsilon q})=O(t^{-\beta})$ for $q$ sufficiently large. 

If $\beta>1$, the fourth term is $O(t^{d-p})=O(t^{-\beta})$
for $p>\beta$ and $d$ sufficiently small.
The case $\beta=1$ differs only by a logarithmic factor so the same choices
of $p$ and $d$ suffice.
If $\beta<1$, the fourth term, ignoring a logarithmic factor, is 
$O(t^{(1-\beta+d)(p+2)-p})=O(t^{-\beta})$
for $p>(2-\beta)/\beta$ and $d$ small.
\end{proof}

Theorem~\ref{thm-exp} is immediate from Corollary~\ref{thm-eigen}
since it is known that the existence of
approximate eigenfunctions implies the periodic data criterion~\eqref{eq-data},
see~\cite[Theorem~1.8]{FMTapp}.

\subsection{Estimates for the Gibbs-Markov map $F:Y\to Y$}

Let $R$ denote the transfer operator for the Gibbs-Markov map $F:Y\to Y$.
We continue to let $f':\Delta'\to\Delta'$ denote the
truncation of $f:\Delta\to\Delta$ with return time $r'=\min\{r,N\}$.
Note that the return map $F=f^r=(f')^{r'}:Y\to Y$ is independent of $N$
and so the operator $R$ is fixed throughout.  
Define $H'(y)=\sum_{\ell=0}^{r'(y)-1}h\circ f^j(y)$.

For $s,z\in\C$, define the twisted transfer operator $R_{s,z}$ to be
$R_{s,z}v=R(e^{sH'}e^{zr'}v)$.

\begin{prop} \label{prop-C6}
$\sum_{j\ge1}|1_{Y_j}H'|_\theta\mu_Y(Y_j)\le |h|_\theta\,\bar r$.
\end{prop}

\begin{proof}
This follows from the estimate
$|1_{Y_j}H'|_\theta\le (r'|Y_j)\,|h|_\theta\le (r|Y_j)|h|_\theta$.
\end{proof}

It follows from Proposition~\ref{prop-C6} that
the estimates in~\cite[Proposition~3.7]{Mapp}
hold independent of $N$.   In particular, 
\begin{align} \label{eq-C6}
|R_{ib,i\omega}^nv|_\theta\le C\{|b||v|_\infty+\theta^n|v|_\theta\},
\end{align}
for all $n,N\ge1$, $|b|>1$, $\omega\in[0,2\pi)$, $v\in F_\theta(Y)$.

Define the norm 
$\|v\|_b=\max\{|v|_\infty,|v|_\theta/(2C|b|)\}$
where $C$ is the constant in~\eqref{eq-C6}.

\begin{lemma} \label{lem-b}
(cf.~\cite[Lemma~3.5]{Mapp}).
Assume no approximate eigenfunctions on $Z_0$ and choose $N$ sufficiently
large that $r|Z_0\le N$.  
Then there exist $\alpha>0$, $C\ge1$ independent of $N$ such that
\begin{align*}
\|(I-R_{ib,i\omega})^{-1}\|_b \le C|b|^\alpha,
\end{align*}
for all $|b|>1$, $\omega\in[0,2\pi)$.
\end{lemma}

\begin{proof}
Since $r|_{Z_0}\le N$,
it makes no difference whether we define $M_{b,\omega}$ using $H$ or $H'$
for the assumption that there are no approximate eigenfunctions on $Z_0$.

When $\omega=0$, it remains to verify that
the proof of~\cite[Lemmas~3.12 and~3,13]{Mapp} goes through unchanged.
The main issue is the dependence on the constant called $C_6$ in~\cite{Mapp}
which potentially
depends on $N$.   However, this constant is shown to be uniform
in~\eqref{eq-C6}.
The remaining arguments in~\cite{Mapp} indeed go through without change
proving the result for $\omega=0$.

As in~\cite[Section~3.3]{Mapp}, the case $\omega\neq0$ presents no
additional complications.
\end{proof}

\begin{prop}   \label{prop-s}
(cf.~\cite[Proposition~3.10]{Mapp}).   
\begin{align*}
\|R_{s,z}-R_{ib,i\omega}\|_b\le Cd_N(|a|+|\sigma|)e^{(|a||h|_\infty+|\sigma|)N},
\end{align*}
for all $s=a+ib$, $z=\sigma+i\omega\in\C$.
\end{prop}

\begin{proof}
The key estimate is~\cite[Proposition~3.9(d)]{Mapp} which states that
\[
\|(R_s-R_{ib})1_{Y_j}\|_b\le C|a|\|1_{Y_j}H'\|_\theta(1+|1_{Y_j}H'|_\theta)
e^{|a||1_{Y_j}H'|_\infty} \mu_Y(Y_j).
\]
Similarly,
\[
\|(R_{s,z}-R_{ib,i\omega})1_{Y_j}\|_b\le 
C(|a|\|1_{Y_j}H'\|_\theta+|\sigma|r'(j))
(1+|1_{Y_j}H'|_\theta)
e^{|a||1_{Y_j}H'|_\infty}e^{|\sigma|r'(j)} \mu_Y(Y_j).
\]
It follows that for each $j\ge1$,
\[
\|(R_{s,z}-R_{ib,i\omega})1_{Y_j}\|_b\le 
C(|a|\|h\|_\theta+|\sigma|)(1+|h|_\theta)
e^{(|a||h|_\infty+|\sigma|)N}\,r'(j)^2\mu_Y(Y_j).
\]
Now $\sum_{j\ge1}r'(j)^2\mu_Y(Y_j)=\sum_{k=1}^Nk^2\mu_Y(r'=k) \le 2d_N$,
so the result follows from the fact that
$R_{s,z}-R_{ib,i\omega}=\sum_{j\ge1}(R_{s,z}-R_{ib,i\omega})1_{Y_j}$.
\end{proof}

In the sequel, $s$ always denotes $s=a+ib\in\C$, similarly
$z=\sigma+i\omega\in\C$.  All constants $C$, $\epsilon$, etc are uniform in 
$|b|>1$ and 
$\omega\in[0,2\pi)$ but we suppress the domain of $b$ and $\omega$.

\begin{lemma} \label{lem-R}
(cf.~\cite[Lemma~3.14]{Mapp}).
Assume no approximate eigenfunctions on $Z_0$ and choose $N$ sufficiently
large that $r|Z_0\le N$.    Let $d>0$ and set $\widetilde{d}_N=d_N N^d$.
There exist $\alpha>0$, $\epsilon>0$ and $C\ge1$ independent of $N$,
such that
\begin{align} \label{eq-R}
\|(I-R_{s,z})^{-1}\|_b \le C|b|^\alpha,
\end{align}
for all $a$, $\sigma\in U_b$, where
\begin{align*}
U_b &= \{a\in\R: |a|<\epsilon\min\{N^{-1}\ln N,\widetilde{d}_N^{\,-1} |b|^{-\alpha}\}\}.
\end{align*}
\end{lemma}

\begin{proof}
Choose $\epsilon<d(|h|_\infty+1)^{-1}$.
It follows from Proposition~\ref{prop-s} that for $s,z$ in the stipulated 
region,
\[
\||R_{s,z}-R_{ib,i\omega}\|_b\le C\epsilon N^{-d}
|b|^{-\alpha} e^{(|h|_\infty+1)\epsilon\ln N}\le C\epsilon |b|^{-\alpha}.
\]
By Lemma~\ref{lem-b},
$\|R_{s,z}-R_{ib,i\omega}\|_b\|(I-R_{ib,i\omega})^{-1}\|_b\le\frac12$ say for
$\epsilon$ sufficiently small.
Using a resolvent inequality as in~\cite[Section~2]{Dolgopyat98b}, we obtain
$\|(I-R_{s,z})^{-1}\|_b\le 2\|(I-R_{ib,i\omega})^{-1}\|_b$ giving the
required result.
\end{proof}

\subsection{Operator renewal sequences}
\label{sec-renewal}

Let $L$ denote the transfer operator for the truncated tower map
$f':\Delta'\to\Delta'$.
Recall that for $s\in\C$, the twisted transfer operator $L_s$ is defined
to be $L_sv=L(e^{sh'}v)$.  Hence
$(L_s^nv)(x)=\sum_{(f')^nz=x}g_n'(z)e^{sh_n'(z)}v(z)$
where $h_n'(z)=h(z)+h(f'z)+\dots+h((f')^{n-1}z)$
and $g_n'(z)$ is the inverse of the Jacobian of $(f')^n$ at $z$.

Let $Z_n=\{y\in Y:r'=n\}$.   Then $\{Z_1,\dots,Z_N\}$
is a finite partition of $Y$.

For $s\in\C$, define the operator renewal sequences 
\[
T_{s,n}=1_YL_{s}^n1_Y, \quad 
R_{s,n}=1_YL_{s}^n1_{Z_n},
\]
and the Fourier series $T,R:\C\to L(F_\theta(Y))$ given by
\[
T_s(z)=\sum_{n=0}^\infty T_{s,n}e^{zn}, \quad
R_s(z)=\sum_{n=1}^N R_{s,n}e^{zn}.
\]
We have the renewal equation $T_s(z)=(I-R_s(z))^{-1}$.
Note also that $R_s(z)v=R_{s,z}v=R(e^{sH'}e^{zr'}v)$.

\begin{lemma} \label{lem-T}
(cf.~\cite[Lemma~4.3]{Mapp})
Assume no approximate eigenfunctions on $Z_0$ and choose $N$ sufficiently
large that $r|Z_0\le N$.
There exist constants $\epsilon,\delta>0,\alpha>0$, $C\ge1$ independent
of $N$ such that
\begin{align} \label{eq-T}
\|T_{s,n}\|_b\le C|b|^\alpha e^{-n\delta\min\{N^{-1}\ln N,\,\widetilde{d}_N^{-1}|b|^{-\alpha}\}},
\end{align}
for all $n\ge1$, and $\sigma\in U_b$.
\end{lemma}

\begin{proof}
By the renewal equation and~\eqref{eq-R},
\begin{align*} 
\|T_s(z)\|_b=\|(I-R_s(z))^{-1}\|_b \le C|b|^\alpha,
\end{align*}
for $a$, $\sigma\in U_b$.

By definition, $R_s(z)$ is a polynomial of degree $N$ in $e^z$
and hence analytic in $z$.   It follows that $T_s(z)$ is analytic in $z$
on the domain of $(I-R_s(z))^{-1}$, namely $U_b$.
Hence the Fourier coefficients $T_{s,n}$ decay at the required rate for
any $\delta<\epsilon$.
\end{proof}

\begin{lemma} \label{lem-L}
(cf.~\cite[Lemma~4.4]{Mapp})
Assume no approximate eigenfunctions on $Z_0$ and choose $N$ sufficiently
large that $r|Z_0\le N$.
Let $d>0$ and set $\widetilde{d}_N=d_N N^d$.
There exist constants $\epsilon,\delta>0,\alpha>0$, $C\ge1$ independent
of $N$ such that
\begin{align} \label{eq-sum}
\sum_{n\ge1}|L_s^nv|_1\le C\|v\|_b\,\widetilde{d}_N 
|b|^\alpha\max\{N(\ln N)^{-1},\widetilde{d}_N |b|^\alpha\},
\end{align}
for all $v\in F_\theta(\Delta')$, $n\ge1$, and $a\in U_b$.
\end{lemma}

\begin{proof}
Recall that $(L_s^nv)(x)=\sum_{(f')^nz=x}g_n'(z)e^{sh_n'(z)}v(z)$.
Following the proof and notation of~\cite[Lemma~4.4]{Mapp}, we write
\[
\SMALL L_s^n=\sum_{i+j+k=n}A_{s,i}T_{s,j}B_{s,k} \enspace+\enspace E_{s,n},
\]
where 
\begin{alignat*}{3}
 (T_{s,n}v)(x)& =\sum_{\substack{f^nz=x \\ x,z\in Y}}
 & \qquad &, &\qquad (A_{s,n}v)(x) &=\sum_{\substack{ f^nz=x \\ z\in Y;\;fz\not\in Y,\ldots,f^nz\not\in Y }}\qquad, \\
  (E_{s,n}v)(x) & =\sum_{\substack{ f^nz=x \\ z\not\in Y,\ldots,f^nz\not\in Y }}
 & \qquad &,& \qquad (B_{s,n}v)(x)& =\sum_{\substack{ f^nz=x \\ z\not\in Y,\ldots,f^{n-1}z\not\in Y;\;
f^nz\in Y }} \qquad,
\end{alignat*}
and we have suppressed the summands $g_n'(z)e^{sh_n'(z)}v(z)$.
We view these as operators $L_s^n:F_\theta(\Delta')\to L^1(\Delta')$,
$T_{s,n}:F_\theta(Y)\to L^\infty(Y)$,
$A_{s,n}:L^\infty(Y)\to L^1(\Delta')$,
$B_{s,n}:F_\theta(\Delta')\to F_\theta(Y)$,
$E_{s,n}:F_\theta(\Delta')\to L^1(\Delta')$,
with the $\|\,\|_b$ norm on $F_\theta(Y)$ and $F_\theta(\Delta')$.
In the corresponding operator norms, we have
\[
\|L_s^n\|\le\sum_{i+j+k=n}\|A_{s,i}\|\,\|T_{s,j}\|\,\|B_{s,k}\|
\enspace+\enspace\|E_{s,n}\|.
\]

Due to truncation, 
$A_{s,n}=B_{s,n}=E_{s,n}=0$ for $n>N$.
We claim further that
\begin{align} \label{eq-ABE}
\|A_{s,n}\| & \le CN^{\epsilon'}\mu_Y(r'\ge n),\quad 
\|B_{s,n}\| \le CN^{\epsilon'} n\mu_Y(r'\ge n),\\
\|E_{s,n}\| & \le \SMALL CN^{\epsilon'}\sum_{k=n}^N\mu_Y(r'\ge k), \nonumber
\end{align}
where $\epsilon'=\epsilon|h|_\infty$,
for all $n\ge 1$, $|a|\le\epsilon N^{-1}\ln N$, $b\in\R$.

Let $u_n=\mu_Y(r'\ge n)$ and $v_n=e^{-cn}$ where $0<c\le \delta N^{-1}\ln N$.
Then 
$(u*v)_n= \sum_{k=1}^N \mu_Y(r'\ge k)e^{-c(n-k)}
\le N^\delta e^{-cn} \sum_{k=1}^N \mu_Y(r'\ge k) 
\le N^\delta\bar r e^{-cn}$.
Similarly, if $u_n'=n\mu_Y(r'\ge k)$, then
$(u'*v)_n\le N^\delta e^{-cn}\sum_{k=1}^N k\mu_Y(r'\ge k)
\le N^\delta d_N e^{-cn}$.
Using this calculation and estimates~\eqref{eq-T},~\eqref{eq-ABE}, we obtain 
\begin{align} \label{eq-L}
\|L_s^n\|\le Cd_N N^{2\epsilon'+2\delta}
|b|^\alpha e^{-n\delta\min\{N^{-1}\ln N,\,\widetilde{d}_N^{-1}|b|^{-\alpha}\}}
+ \|E_{s,n}\|.
\end{align} 
Moreover, $\sum_{n\ge 1}\|E_{s,n}\|
\le CN^\epsilon\sum_{n=1}^N\sum_{k=n}^N\mu_Y(r'\ge k)=CN^\epsilon d_N$.
Shrink $\epsilon$ and $\delta$ if necessary so that $2\epsilon'+2\delta<d$.
Since $(1-e^{-x})^{-1}\le 2x^{-1}$ for $x>0$ small,
we obtain the required estimate for $\sum_{n\ge 1}L_s^n$,

It remains to verify estimates~\eqref{eq-ABE}.  Note that the support of
$A_{s,n}v$ is contained in level $n\le N$ of the tower and has measure at 
most $\sum_{r'(j)>n}\mu_{\Delta'}(Y_j)\le (1/\bar r')\mu_Y(r'\ge n)$.  
For $x$ in level $n$, we have $(A_{s,n}v)(x)=e^{sh'_n(z)}v(z)$
where $z$ is the unique point in $Y$ with $(f')^nz=x$, and so
$|A_{s,n}v|_\infty\le e^{\epsilon N^{-1}\ln N\,n|h|_\infty}|v|_\infty
\le N^{\epsilon|h|_\infty}|v|_\infty$.
Hence $|A_{s,n}v|_1\le (1/\bar r')N^{\epsilon'}\mu_Y(r'\ge n)$.  Similarly,
\[
|E_{s,n}v|_1\le N^{\epsilon'}\sum_{\substack{ r'(j)>n\\n< \ell<r'(j)}}
\mu_{\Delta'}(\Delta_{j,\ell})|v|_\infty \le (1/\bar r')N^{\epsilon'}
\sum_{k=n+2}^N\mu_Y(r'\ge k)\|v\|_b.
\]
Finally, if $y\in Y$, then
$(B_{s,n}v)(y)=\sum_{r'(j)>n}g_n'(z_j')e^{sh'_n(z_j')}v(z_j')$
where $z_j'$ is the unique preimage of $y$ in $\Delta_{j,r'(j)-n}$.
Since $f':\Delta_{j,\ell}\to\Delta_{j,\ell+1}$ is an isometry 
for $\ell<r'-1$, we can write $g_n'(z_j')=g(z_j)$ where $z_j$ is the 
unique point satisfying $z_j\in Y_j$, $Fz_j=y$,
and $g$ is the Jacobian in part (4) of the definition of nonuniform
expansion in Section~\ref{sec-main-a}.
(Alternatively, $g$ is the weight in the definition
$(Rv)(x)=\sum_{Fy=x}g(y)v(y)$ of the transfer operator $R$ for the 
Gibbs-Markov map $F:Y\to Y$.)
The log-H\"older condition on $g$ implies that $|g(y)|\le C\mu_Y(Y_j)$
and $|g(y)/g(\hat y)-1|\le Cd_\theta(y,\hat y)$ for all $y,\hat y\in Y_j$.
Hence $|B_{s,n}v|_\infty\le \sum_{r'(j)>n}C\mu_Y(Y_j)N^{\epsilon'}|v|_\infty
\le CN^{\epsilon'}\mu_Y(r'\ge n)|v|_\infty$ and 
$|(B_{s,n}v)(y)-B_{s,n}v)(\hat y)|
\le\sum_{r'(y)>n}|g(z_j)e^{sh_n'(z_j')}v(z_j')-
g(\hat z_j)e^{sh_n'(\hat z_j')}v(\hat z_j')|
\le d_\theta(y,\hat y)\sum_{r'(y)>n}
\Bigl(C\mu_Y(Y_j)N^{\epsilon'}|v|_\theta
+C\mu_Y(Y_j)N^{\epsilon'}|s|n|h|_\theta|v|_\infty
+C\mu_Y(Y_j)N^{\epsilon'}|v|_\infty\Bigr)$
so that
$|B_{s,n}v|_\theta\le CN^{\epsilon'} |b|n\mu_Y(r'\ge n)|v\|_b$.
It follows that 
$\|B_{s,n}\|_b\le CN^{\epsilon'}n\mu_Y(r'\ge n)$ completing the verification
of estimates~\eqref{eq-ABE}.
\end{proof}

\begin{pfof}{Theorem~\ref{thm-eigen}}
By the formula in Step~4, Section~\ref{sec-main-d}, it follows from
Lemma~\ref{lem-L} that
\begin{align} \label{eq-hatrho}
\hat\rho(s)\le C\|v\|_\theta|w|_\infty \widetilde{d}_N 
|b|^\alpha\max\{N(\ln N)^{-1},\widetilde{d}_N |b|^\alpha\}.
\end{align}
To recover $\rho'(t)$ we integrate along a contour $a=a_N(b)$ in the 
left-half-plane; specifically 
$a=\epsilon\min\{N^{-1}\ln N,\widetilde{d}_N^{\,-1}b^{-\alpha}\}$ for $b>1$
(and decreased $\epsilon$).  
Integrating by parts $m$ times as in~\cite{Dolgopyat98b}, we 
obtain the integrals
\[
I=\int_1^\infty \widetilde{d}_NNe^{-\epsilon N^{-1}\ln N\, t}
b^{\alpha-m}\,db, \quad 
II=\int_1^\infty e^{-\epsilon \widetilde{d}_N^{\,-1}b^{-\alpha}t}
\widetilde{d}_N^{\,2} b^{2\alpha-m}\,db.
\]
Taking $m>\alpha+1$ yields 
$I\le C\widetilde{d}_NNe^{-\epsilon N^{-1}\ln N\, t}$.
A change of variables yields
\begin{align*}
II & \le \alpha^{-1}\epsilon^{-p}t^{-p} \widetilde{d}_N^{\,p+2}\int_0^\infty e^{-y} y^{p-1}\,dy 
=\alpha^{-1}\epsilon^{-p}(p-1)!\, \widetilde{d}_N^{\,p+2} t^{-p},
\end{align*}
with $m=(p+2)\alpha+1$.
\end{pfof}

\section{Decay for nonuniformly hyperbolic flows}
\label{sec-hyp}

In this section we prove Theorem~\ref{thm-hyp}.
The main steps are the same as for Theorem~\ref{thm-exp}, but
there is an additional step 
between Steps 3 and 4 where we pass from the Young tower
to a nonuniformly expanding quotient tower 
$\bar{f'}:\bar{\Delta'}\to\bar{\Delta'}$
by quotienting along stable manifolds.

In Subsection~\ref{sec-hyp-a}, we include the necessary background material and
notation from Young~\cite{Young98,Young99} on nonuniformly hyperbolic 
diffeomorphisms and towers.
In Subsection~\ref{sec-hyp-b}, we use approximation arguments to reduce
the nonuniformly hyperbolic case to the nonuniformly expanding case
studied in Section~\ref{sec-exp}.

\subsection{Background on nonuniformly hyperbolic systems}
\label{sec-hyp-a}

Let $T:M\to M$ be a nonuniformly hyperbolic diffeomorphism in the sense of
Young~\cite{Young98,Young99}. As described in Section~\ref{sec-main-b},
there is a partition $\{Y_j\}$ of $Y\subset M$
with return time function $r:Y\to\Z^+$, constant on partition
elements $\{Y_j\}$, and induced return map
$F:Y\to Y$ given by $F(y)=T^{r(y)}(y)$.  There exists an
ergodic $T$-invariant probability measure $\nu$ that is an SRB measure.

Let $\Delta=\{(y,\ell):y\in Y,\,\ell=0,\dots,r(y)-1\}$ and
define the tower map $f:\Delta\to\Delta$ by setting $f(y,\ell)=(y,\ell+1)$
for $0\le \ell<r(y)-1$ and $f(y,r(y)-1)=(Fy,0)$.
The projection $\pi:\Delta\to M$ given by $\pi(y,\ell)=T^\ell y$
is a semiconjugacy between $f:\Delta\to\Delta$ and $T:M\to M$.

The subset $Y$ is covered by families of {\em stable disks}
$\{W^s(y),\,y\in Y\}$
and {\em unstable disks} $\{W^u(y),\,y\in Y\}$ such that
each stable disk intersects each unstable disk in exactly one point.
For $p=(x,\ell),\,q=(y,\ell)\in \Delta$,
we write $q\in W^s(p)$ if $y\in W^s(x)$ (and
$q\in W^u(p)$ if $y\in W^u(x)$).

Quotienting out the stable directions, we obtain the quotient
maps $\bar f:\bar\Delta\to\bar\Delta$ and $\bar F:\bar Y\to\bar Y$.

\begin{prop}[ {\cite{Young98,Young99} } ]
The quotient tower map $\bar f:\bar\Delta\to\bar\Delta$
is a nonuniformly expanding tower map of the type
considered in Section~\ref{sec-exp}.
In particular, there are $\bar F$ and $\bar f$-invariant
measures $\bar\mu$ and $\bar\mu\times\mu_C/\int_{\bar Y}r\,d\bar\mu$
on $\bar Y$ and $\bar\Delta$ respectively, such that
$\bar F:\bar Y\to\bar Y$ is Gibbs-Markov with respect to the
quotient partition $\{\bar Y_j\}$.
Moreover, there is a $f$-invariant measure $\mu$ on $\Delta$ such that
the natural projection $\bar\pi:\Delta\to\bar\Delta$ and the projection
$\pi:\Delta\to M$ are measure-preserving semiconjugacies.
\qed
\end{prop}

In Step~1, we defined 
$s:\bar\Delta\times\bar\Delta\to\N$ relative to returns under $\bar F$
to the partition $\{\bar Y_j\}$.
This lifts to a separation time 
$s:\Delta\times\Delta\to\N$ given by $s(p,q)=s(\bar\pi p,\bar\pi q)$.
Note that $s$ is defined on both $\bar\Delta$
and $\Delta$, but the metric $d_\theta(p,q)=\theta^{s(p,q)}$ is defined only 
on $\bar\Delta$.

We assume that there exists $\gamma\in(0,1)$ such that
\begin{itemize}
\item[(P1)]  If $q\in W^s(p)$, then
$d(\pi f^np,\pi f^nq)\le C\gamma^n$ for all $n\ge1$.
\item[(P2)]  If $q\in W^u(p)$, then
$d(\pi f^np,\pi f^nq)\le C\gamma^{s(p,q)-n}$ for $0\le n< s(p,q)$.
\end{itemize}
This means that there is exponential contraction along stable disks but
nonuniform expansion along unstable disks.  

\begin{prop} \label{prop-W}
$d(T^n\pi p,T^n\pi q)\le C\gamma^{\min\{n,s(p,q)-n\}}$ for all
$p,q\in\Delta$, $0\le n\le s(p,q)$.
\end{prop}

\begin{proof}
Define $z=W^s(p)\cap W^s(q)$.
By (P1),
$d(\pi f^np,\pi f^nz)\le C\gamma^n$.
Moreover, $s(z,q)=s(p,q)$ and so by (P2),
$d(\pi f^nz,\pi f^nq)\le C\gamma^{s(p,q)-n}$.
\end{proof}

\subsection{Proof of Theorem~\ref{thm-hyp}}
\label{sec-hyp-b}

We continue to assume that $T:M\to M$ is a nonuniformly hyperbolic 
diffeomorphism,
modelled by a Young tower $f:\Delta\to\Delta$ as in Subsection~\ref{sec-hyp-a}.
We have the measure-preserving semiconjugacy $\pi:\Delta\to M$.

Let $h:M\to\R^+$  be a $\eta$-H\"older roof function with
associated suspension flow $T_t:M^h\to M^h$.
Define $\tilde h=h\circ \pi$ with suspension flow
$f_t:\Delta^{\tilde h}\to\Delta^{\tilde h}$.  The projection
$\pi:\Delta^{\tilde h}\to M^{h}$ defined by $\pi(p,u)=(\pi p,u)$
is a measure-preserving semiconjugacy.
Given $v,w\in C^\eta(M^h)$,
let $\tilde v=v\circ \pi$, $\tilde w=w\circ \pi$.
It suffices to prove decay of correlations for the
observations $\tilde v,\tilde w :\Delta^{\tilde h}\to\R$.

We now introduce the truncation $r'=\min\{r,N\}$ and the truncated
tower map $f':\Delta'\to\Delta'$.   
The argument in Section~\ref{sec-trunc} for the nonuniformly expanding
case applies equally to the the nonuniformly hyperbolic case.  

To simplify notation, in the remainder of this
section we write $f:\Delta\to\Delta$ for the truncated tower map
and $\mu$ for the measure on $\Delta'$.
Note that estimate (P2) is unaffected by truncation since the
return map to $Y$ is unchanged.   Also, estimate (P1) can only be improved
by truncation.
In particular, Proposition~\ref{prop-W} remains valid.
We note that many of the objects defined below, such as $\chi$,
$v_{s,k}$ and so on, depend on $N$.  However, the estimates involve
universal constants independent of $N$.

As in the uniformly expanding case, the significant part of the Laplace 
transform of the correlation function for the truncated flow has the form
\[
\hat\rho(s)=\sum_{n\ge1} \int_\Delta e^{-s\tilde h_n}v_s
\; w_s \circ f^n \,d\mu,
\]
where $v_s(p)=\int_0^{\tilde h(p)}e^{su}\tilde v(p,u)du$ and
$w_s(p)=\int_0^{\tilde h(p)}e^{-su}\tilde w(p,u)du$.

To estimate $\hat\rho(s)$, 
the first step is to write $\tilde h$ as a coboundary
plus a roof function that ``depends only on future coordinates''.

\begin{lemma} \label{lem-future}
There exist functions $\bar h,\chi:\Delta\to\R$ such that
\begin{itemize}
\item[(i)] $\tilde h=\bar h+\chi-\chi\circ f$,
\item[(ii)] $\chi\in L^\infty(\Delta)$ and $|\chi|_\infty\le C$
(independent of $N$),
\item[(iii)] If $s(p,q)\ge3k$, then $|\chi(f^kp)-\chi(f^kq)|\le C\gamma_1^k$, where $\gamma_1=\gamma^\eta$,
\item[(iv)] $\bar h(p)=\bar h(q)$ for all $p\in W^s(q)$,
\item[(v)] $\bar h:\bar\Delta\to\R$ is Lipschitz with respect to the
metric $d_\theta$, for $\theta=\gamma_1^{1/3}$.
\end{itemize}
\end{lemma}

\begin{proof}
We modify the proof of~\cite[Lemma~5.4]{Mapp}.
Instead of the two separation times $s$ and $s_1$ in~\cite{Mapp}, we have
only the separation time $s$.   Most of~\cite[Lemma~5.4]{Mapp} goes through 
word for word with $s$ substituted for $s_1$.
The proof only differs in part (e): instead of choosing $p,q\in\Delta$
with $s_1(p,q)\ge 2k+1$, we require that $s(p,q)\ge 3k+1$.
Since $\bar\Delta$ is Markov, we can
choose $\bar p'\in \bar f^{-k}\bar p$, $\bar q'\in \bar f^{-k}\bar q$ 
with $s(p',q')\ge 3k+1$.
(Unlike in~\cite{Mapp}, the separation of $p,q$ does not necessarily
increase with each backward iterate.)
By~(i),~(iii) and the H\"older continuity of $h$, we have that
$|\bar h(p)-\bar h(q)|=|\bar h(f^kp')-\bar h(f^kq')|
\le C\gamma_1^k$ as required.
\end{proof}

\noindent
(The proof of Lemma~\ref{lem-future} shows that the introduction 
of $s_1$ in~\cite{Mapp} is unnecessary.)

By Lemma~\ref{lem-future}, we can write
$\hat\rho(s)=\sum_{n\ge1} \int_\Delta e^{-s\bar h_n}(e^{-s\chi}v_s)
\; (e^{s\chi}w_s) \circ f^n \,d\mu$.
Next we approximate $e^{-s\chi}v_s$ and $e^{s\chi}w_s$
by functions that ``depend only on finitely many coordinates''.
For $k\ge1$, define $v_{s,k}(p)=\inf\{(e^{-s\chi}v_s)(f^kq):s(p,q)\ge3k\}$.

\begin{lemma} \label{lem-vsk}
The function $v_{s,k}:\Delta\to\R$ lies in $L^\infty(\Delta)$ and
projects down to a Lipschitz observation
$\overline{v_{s,k}}:\bar\Delta\to\R$.  
Within the region $s=a+ib$, $|a|\le1$, $|b|\ge1$,
\begin{itemize}
\item[(a)] $|\overline{v_{s,k}}|_\infty = |v_{s,k}|_\infty \le
e^{|\chi|_\infty}|v_s|_\infty \le C|\tilde v|_\infty = C|v|_\infty$.
\item[(b)] $|\overline{v_{s,k}}|_\theta  \le
C|v|_\infty\,\theta^{-3k}$.
\item[(c)]  $|(e^{-s\chi}v_s)\circ f^k - v_{s,k}|_\infty  \le
C\|v\|_{\eta}|b|\gamma_1^k$.
\end{itemize}
\end{lemma}

\begin{proof}
The proof is unchanged from~\cite[Lemma~5.5]{Mapp} except that $s$ is 
again substituted for $s_1$.
\end{proof}

Write
$\int_\Delta e^{-s\bar h_n}(e^{-s\chi}v_s)\,(e^{s\chi}w_s)\circ f^n\,d\mu    =
\int_\Delta e^{-s\bar h_n\circ f^k}(e^{-s\chi}v_s)\circ f^k\,(e^{s\chi}w_s)\circ f^k\circ f^n\,d\mu = I_1+I_2+I_3$, where
\begin{align*}
I_1 & = {\SMALL\int}_\Delta e^{-s\bar h_n\circ f^k}(e^{-s\chi}v_s)\circ f^k\,((e^{s\chi}w_s)\circ f^k-w_{s,k})\circ f^n\,d\mu, \\
I_2 & = {\SMALL\int}_\Delta e^{-s\bar h_n\circ f^k}((e^{-s\chi}v_s)\circ f^k-v_{s,k})\,w_{s,k}\circ f^n\,d\mu, \\
I_3 & = {\SMALL\int}_\Delta e^{-s\bar h_n\circ f^k}v_{s,k}\,w_{s,k}\circ f^n\,d\mu.
\end{align*}
By Lemma~\ref{lem-vsk},
\[
|I_1| \le e^{|a||\bar h_n|_\infty}e^{|\chi|_\infty}|v_s|_\infty 
|(e^{s\chi}w_s)\circ f^k-w_{s,k}|_\infty \le
C|v|_\infty \|w\|_{\eta} |b| e^{n|a||h|_\infty}\gamma_1^k,
\]
and similarly $|I_2|\le
C\|v\|_{\eta}|w|_\infty|b| e^{n|a||h|_\infty} \gamma_1^k$.
Hence
\[
|I_1|,|I_2| \le
C\|v\|_{\eta}\|w\|_{\eta} |b| 
e^{n\epsilon\min\{N^{-1}\ln N,\,\widetilde{d}_N^{\,-1} |b|^{-\alpha}\}|h|_\infty}
\gamma_1^k,
\]
for all $a\in U_b$.

The integrand in $I_3$ projects down to $\bar\Delta$ and
$\bar h_n\circ \bar f^k=\bar h_n+\bar h_k\circ \bar f^n-\bar h_k$, so 
\[
I_3={\SMALL\int}_{\bar \Delta} e^{-s\bar h_n}[e^{s\bar h_k}\overline{v_{s,k}}]\;[e^{-s\bar h_k}\overline{w_{s,k}}]\circ \bar f^n\,d\bar\mu
={\SMALL\int}_{\bar \Delta} L_{-s}^n[e^{s\bar h_k}\overline{v_{s,k}}]\;
[e^{-s\bar h_k}\overline{w_{s,k}}]\,d\bar\mu.
\]
Here, $L$ is the transfer operator for the truncated quotient tower map
$\bar f:\bar\Delta\to\bar\Delta$, and $L_su= L(e^{s\bar h}u)$.  By~\eqref{eq-L},
\[
\|L_s^n\|\le C\widetilde{d}_N 
|b|^\alpha e^{-n\delta\min\{N^{-1}\ln N,\,\widetilde{d}_N^{\,-1}|b|^{-\alpha}\}}
+\|E_{s,n}\|
\]
on $F_\theta(\bar\Delta)$.  Hence,
\begin{align*}
|I_3| &\le \|L_{-s}^n\|_b
\|e^{s\bar h_k}\|_\theta\|\overline{v_{s,k}}\|_\theta|e^{-sh_k}|_\infty|\overline{w_{s,k}}|_\infty \\
&\le
 C|v|_\infty|w|_\infty (\widetilde{d}_N|b|^{\alpha} 
e^{-n\delta\min\{N^{-1}\ln N,\,\widetilde{d}_N^{\,-1}|b|^{-\alpha}\}}+\|E_{s,n}\|)
|b|\theta^{-4k}e^{2k|h|_\infty}.
\end{align*}

Choose $k=k(b,n,N)$ so that
\[
(e^{2|h|_\infty}\theta^{-4})^k\sim 
e^{\frac 12 n\delta\min\{N^{-1}\ln N,\,\widetilde{d}_N^{\,-1}|b|^{-\alpha}\}}.
\]
Then there exists $\delta'>0$ (depending on $\gamma_1$ and $\theta$)
such that
\begin{align*}
& I_1,I_2 =O(e^{-n(\delta'-\epsilon)\min\{N^{-1}\ln N,\,\widetilde{d}_N^{\,-1}|b|^{-\alpha}\}|h|_\infty}|b|), \\
& I_3 =O(\widetilde{d}_Ne^{-\frac12 n\delta\min\{N^{-1}\ln N,\,\widetilde{d}_N^{\,-1} |b|^{-\alpha}\}}
|b|^{\alpha+1}) +O(N^{\delta/2}\|E_{s,n}\||b|).
\end{align*}
Here, we have used the fact that $E_{s,n}=0$ for $n>N$.
Choosing $\epsilon$ small enough, we obtain a new $\delta>0$ such that
\begin{align*}
& |{\SMALL\int}_\Delta e^{-s\bar h_n}(e^{-s\chi}v_s) \; (e^{s\chi}w_s) \circ f^n d\mu| \\ & \qquad  \le
C\|v\|_{\eta}\|w\|_{\eta}
(\widetilde{d}_N e^{-n\delta\min\{N^{-1}\ln N,\,\widetilde{d}_N^{\,-1}|b|^{-\alpha}\}}+N^\delta\|E_{s,n}\|)
|b|^{\alpha+1}.
\end{align*}
Summing over $n$ as in Lemma~\ref{lem-L}, we obtain
\[
|\hat\rho(s)| \le C\|v\|_{\eta}\|w\|_{\eta}\,\widetilde{d}_N|b|^{\alpha+1}\max\{N(\ln N)^{-1},\widetilde{d}_N |b|^\alpha\}.
\]
This is almost identical to the estimate~\eqref{eq-hatrho} obtained 
in the nonuniform expanding case, and so we recover the
required decay of correlation result for $\rho'(t)$ in Theorem~\ref{thm-eigen}
as before (but with $m=(p+2)\alpha+2$) and hence for $\rho(t)$.
\qed

\section{Decay for flows with unbounded roof functions}
\label{sec-unbounded}

In this section, we prove Theorem~\ref{thm-unbounded}.
The strategy is similar to the previous sections but now two
truncations are required.  Let 
\begin{align*}
U_b &= \bigl\{a\in\R:
|a|<\epsilon\min\{N^{-1},\widetilde{d}_N^{\,-1}|b|^{-\alpha}\}\bigr\}, \\[.75ex]
V_b &=\bigl\{\sigma\in\R:
|\sigma|<\epsilon\widetilde{d}_N^{\,-1}|b|^{-\alpha}\bigr\},
\end{align*}
where $\widetilde{d}_N$ is defined below in Lemma~\ref{lem-RR}.
In the semiflow case, the steps are as follows:
\begin{itemize}
\item[(a)]  Model by a tower.
\item[(b)]  Truncate $h$ to $h'=\min\{h,N\}$.
This leads to an error $O(N^{-\beta}+tN^{-(\beta+1)})$ in the 
correlation function.
\item[(c)]  Truncate $r$ to $r'=\min\{r,[q\ln N]\}$.
The truncation error takes the form $O(tN^{-(cq-1)})$.
Choosing $q>(\beta+2)/c$ ensures that this is dominated by the error in (b).
\item[(d)] $\|(I-R_{s,z})^{-1}\|_b\le C|b|^\alpha$
for $a\in U_b$, $\sigma\in V_b$.
\item[(e)] 
$|L_s^nv|_1\le C\|v\|_bN|b|^\alpha e^{-n\delta\widetilde{d}_N^{-1}|b|^{-\alpha}}$ for $a\in U_b$.
\item[(f)] 
$|\hat\rho(s)|\le C\|v\|_\theta|w|_\infty
 N^3\widetilde{d}_N|b|^{2\alpha}$ for $a\in U_b$,
where $\hat\rho(s)$ is the Laplace transform of the (doubly) truncated 
correlation function $\rho'(t)$.
\item[(g)]  
$|\rho'(t)|\le C\|v\|_{m,\theta}|w|_\infty\{ 
N^3\widetilde{d}_Ne^{-\epsilon N^{-1}t}+ N^3 \widetilde{d}_N^{\,p+1} t^{-p}\}$.
\item[(h)] Specify $N=N(t)$.
\end{itemize}

Step~(a) is identical to
Step~1 in Section~\ref{sec-main-d}. We may assume from now on that 
the nonuniformly expanding map is a tower map $f:\Delta\to\Delta$
and that $h:\Delta\to\R^+$ is a (nonuniformly) Lipschitz roof function.

\subsection{Truncation of $h$}

In this subsection, we carry out Step~(b).  Let 
$\Delta(n)=\bigcup\{\Delta_{j,\ell}:\|h1_{\Delta_{j,\ell}}\|_{\theta}\ge n\}$.
Condition (7) on $h$ in Section~\ref{sec-main-c} guarantees that 
$\mu_\Delta(\Delta(n))\le Cn^{-(\beta+1)}$.

Fix $N\ge1$ and let $h'=\min\{h,N\}$.   We form the suspension flows
$f_t:\Delta^h\to\Delta^h$ and $f'_t:\Delta^{h'}\to\Delta^{h'}$.
Observables $v,w$ on $\Delta^h$ 
restrict to observables on $\Delta^{h'}$ and we define
the correlation functions $\rho(t)$ and $\rho'(t)$.

Write $\Delta^h=\Delta^h_{\rm left}\dot\cup \Delta^h_{\rm right}$ where
\[
\Delta^h_{\rm left}=\{(x,u)\in\Delta^h:h(x)\le N\}, \qquad
\Delta^h_{\rm right}=\{(x,u)\in\Delta^h:h(x)>N\}.
\]
As in Proposition~\ref{prop-r}, we obtain
$\bar h-\bar{h'}\le CN^{-\beta}$ and
$\mu_{\Delta^h}(\Delta^h_{\rm right})\le CN^{-\beta}$.

\begin{prop} \label{prop-Ekk}
For $k\ge1$, define
\[
E_k=\{p\in\Delta^h:f_tp\in\Delta^h_{\rm right} 
\enspace\text{for some $t\in[0,k]$}\}.
\]
Then
$\mu_{\Delta^h}(E_k)\le C\{N^{-\beta}+kN^{-(\beta+1)}\}$
for all $N\ge2$.
\end{prop}

\begin{proof}
Write $E_k$ as the disjoint union $E_k=
\bigcup_{j=1}^k G_j$ where
\[
G_j=\{
f_tp\in\Delta^h_{\rm left}\enspace\text{for $t\in[0,j-1)$}
\enspace\text{and}\enspace
f_tp\in\Delta^h_{\rm right}
\enspace\text{for some $t\in[j-1,j]$}\}.
\]

For $j\ge2$, it follows from the definition that if $p\in G_j$, then
$f_jp\in\Delta^1_{\rm right}$ where
$\Delta^1_{\rm right}=\{(x,u)\in\Delta\times[0,1]:h(x)>N\}$.
Hence $\mu_{\Delta^h}(G_j)\le \mu_{\Delta^h}(f_j^{-1}(\Delta^1_{\rm right}))
=\mu_{\Delta^h}(\Delta^1_{\rm right})=(1/\bar h)\mu_\Delta(h>N)
\le C N^{-(\beta+1)}$.

If $p\in G_1$, then either $p\in\Delta^h_{\rm right}$ or
$f_1 p\in \Delta^1_{\rm right}$.
Hence $\mu_{\Delta^h}(G_1)\le C N^{-\beta}+C N^{-(\beta+1)}$.
\end{proof}

\begin{lemma} 
Suppose that $v,w:\Delta^h\to\R$ lie in $L^\infty$ and define $\rho(t)$, 
$\rho'(t)$ as indicated above.  For $N\ge 2$, $t>0$,
\[
|\rho(t)-\rho'(t)|\le
C|v|_\infty|w|_\infty \{N^{-\beta}+ tN^{-(\beta+1)}\}.
\]
\end{lemma}

\begin{proof}
For notational convenience, we write $\Omega=\Delta^h$ and
$\Omega'=\Delta^{h'}$.
Let $A=\int_\Omega v\,w\circ f_t\,d\mu_\Omega$, $A'=\int_{\Omega'} v\,w\circ f'_t\,d\mu_{\Omega'}$.
Then 
\begin{align*}
A-A' & =\int_\Omega (v\,w\circ f_t-v\,w\circ f'_t)\,d\mu_\Omega
+\Bigl(\int_\Omega v\,w\circ f'_t\,d\mu_\Omega-\int_{\Omega'} v\,w\circ f'_t\,d\mu_{\Omega'}
\Bigr) \\ & = I + II.
\end{align*}
Using Proposition~\ref{prop-Ekk}, we compute that
\[
|I|\le 2|v|_\infty|w|_\infty\mu_\Omega\{f_t\neq f'_t\}
\le C|v|_\infty|w|_\infty\{N^{-\beta}+(t+1)N^{-(\beta+1)}\}.
\]
Next,
\[
II=(1/\bar h)\int_\Delta\int_{h'}^h v\, w\circ f_t'\,du\,d\mu_\Delta
+\Bigl((1/\bar h)-(1/\bar h')\Bigr)\int_\Delta\int_0^{h'} v\, w\circ f_t'\,du\,d\mu_\Delta.
\]
Hence
\begin{align*}
|II|&\le (1/\bar h)|v|_\infty|w|_\infty(\bar h-\bar h')+
(1/\bar h)(1/\bar h')(\bar h-\bar h')|v|_\infty|w|_\infty\bar h' \\
& \le  C |v|_\infty|w|_\infty N^{-\beta}.
\end{align*}
The result follows.
\end{proof}

\subsection{Truncation of $r$}
Recall that $\mu_Y(r>n)=O(e^{-cn})$ where $c>0$.
We make the second truncation $r'=\min\{r,[q\ln N]\}$.
Following Section~\ref{sec-trunc}, we obtain
$\bar r-\bar r'\le C N^{-cq}$, 
$\mu_\Delta(\Delta_{\rm right})\le C N^{-cq}$ and 
$\mu_\Delta(E_k)\le C kN^{-cq}$. 
(Note that none of these calculations depends on $h$.)
The proof of Lemma~\ref{lem-trunc} proceeds as before except that there
is need for care since $|h'|_\infty=N$.   This leads to the loss of one factor
of $N$ and hence the truncation error $O(tN^{-(cq-1)})$.

\subsection{Decay of the truncated correlation function}
\label{sec-semiunb}

In this subsection, we carry out Steps~(d)-(h),
completing the proof of Theorem~\ref{thm-unbounded}.
Let $Y(n)=\bigcup\{Y_j:\|1_{Y_j}H\|_\theta\ge n\}$.

\begin{lemma} \label{lem-Yn}
$\mu_Y(Y(n))=O((\ln n)^{\beta+2}n^{-(\beta+1)})$.
\end{lemma}

\begin{proof}
Let $Q>0$ and write $\mu_Y(Y(n))\le \mu_Y(r>[Q\ln n])+\sum_{k=1}^{[Q\ln n]}
\mu_Y(\{r=k\}\cap Y(n))$.
If $Y_j\subset \{r=k\}\cap Y(n)$, then $\|1_{\Delta_{j,\ell}}h\|_\theta
>n/k$ for some $\ell<k$, and so $\Delta_{j,\ell}\subset\Delta(n/k)$.
Since $\mu_Y(Y_j)=\bar r\mu_\Delta(\Delta_{j,\ell})$,
\[
\mu_Y(\{r=k\}\cap Y(n))\le \bar r\mu_\Delta(\Delta(n/k))
\le C(k/n)^{\beta+1}.
\]
Hence 
\[
\mu_Y(Y(n))\le C \Bigl(e^{-cQ\ln n}+\sum_{k=1}^{[Q\ln n]}(k/n)^{\beta+1} \Bigr)
\le C (e^{-cQ\ln n}+(\ln n)^{\beta+2}n^{-(\beta+1)}).
\]
Now choose $Q=(\beta+1)/c$.
\end{proof}

Since $h'=\min\{h,N\}$, we have $\|1_{\Delta_{j,\ell}}h'\|_\theta\le
\|1_{\Delta_{j,\ell}}h\|_\theta$ for all 
partition elements $\Delta_{j,\ell}$, and hence
$\|1_{Y_j}H'\|_\theta \le \|1_{Y_j}H\|_\theta$.
By Lemma~\ref{lem-Yn},
\[
\sum_{j\ge1}\|1_{Y_j}H'\|_\theta\, \mu_Y(Y_j)\le
\sum_{j\ge1}\|1_{Y_j}H\|_\theta\, \mu_Y(Y_j)<\infty.
\]
This corresponds to Proposition~\ref{prop-C6} and 
guarantees that we obtain a basic inequality
\[
|R_{ib,i\omega}^nv|_\theta\le C\{|b||v|_\infty +\theta^n|v|_\theta\},
\]
uniformly in $N$.
Define $\|v\|_b=\max\{|v|_\infty,|v|_\theta/(2C|b|)\}$.

Let $Z_0\subset Y$ be a finite subsystem and let $Z_0^r$ denote the 
part of the tower $\Delta$ over $Z_0$.   Then $Z_0^r$ 
consists of finitely many partition elements so that $h|Z_0^r$ is bounded.
It follows that there exists $N_0$ such that $h|Z_0^r=h'|Z_0^{r'}$ for all
$N\ge N_0$.  In particular, $H|Z_0=H'|Z_0$ for all $N\ge N_0$.
We now have the ingredients required for the analogue of
Lemma~\ref{lem-b}:

\begin{lemma}   
Assume no approximate eigenfunctions on $Z_0$ and choose 
$N$ sufficiently large that $H|Z_0=H'|Z_0$.
Then there exist $\alpha>0$, $C\ge1$ independent of $N$
such that
\begin{align*} 
\|(I-R_{ib,i\omega})^{-1}\|_b \le C|b|^\alpha,
\end{align*}
for all $|b|>1$, $\omega\in[0,2\pi)$.
\qed
\end{lemma}

\begin{prop}    \label{prop-R}
Let $d_N=\sum_{k=1}^N k\mu_Y(Y(k))$.  Then
\begin{align*} 
\|R_{s,z}-R_{ib,i\omega}\|_b\le Cd_N(|a|+|\sigma|)e^{q(|a|N+|\sigma|)\ln N},
\end{align*}
for all $s$, $z\in\C$.
\end{prop}

\begin{proof}
As in the proof of Proposition~\ref{prop-s}, we have
\[
\|(R_{s,z}-R_{ib,i\omega})1_{Y_j}\|_b\le 
C(|a|\|1_{Y_j}H'\|_\theta+|\sigma|r'(j))
(1+|1_{Y_j}H'|_\theta)
e^{|a||1_{Y_j}H'|_\infty}e^{|\sigma|r'(j)} \mu_Y(Y_j).
\]
Note that $r'(j)\le |\frac1h|_\infty |1_{Y_j}H'|_\infty$.
Since $r_j'\le q\ln N$ and $|1_{Y_j}H'|_\infty\le qN\ln N$,
\begin{align*}
& \SMALL \|(R_{s,z}-R_{ib,i\omega})1_{Y_j}\|_b
\SMALL \le C(|a|+|\sigma|)e^{q(|a|N+|\sigma|)\ln N}
\|1_{Y_j}H'\|^2_\theta \mu_Y(Y_j).
\end{align*}
Now sum over $j\ge1$.
\end{proof}

Combining these two results leads to the following analogue of
Lemma~\ref{lem-R}, completing Step~(d):

\begin{lemma} \label{lem-RR}
Assume no approximate eigenfunctions on $Z_0$ and choose 
$N$ sufficiently large that $H|Z_0=H'|Z_0$.
Let $d>0$ and set $\widetilde{d}_N=d_N N^d$.
Define $U_b$, $V_b$ as at the beginning of the section.
Then there exist $\alpha>0,\epsilon>0$ and $C\ge1$ independent of $N$
such that
\[
\|(I-R_{s,z})^{-1}\|_b \le C|b|^\alpha,
\]
for all $a\in U_b$, $\sigma\in V_b$.
\qed
\end{lemma}

Next we carry out Step~(e).
\begin{lemma} \label{lem-LL}
Assume no approximate eigenfunctions on $Z_0$ and choose 
$N$ sufficiently large that $H|Z_0=H'|Z_0$.
Let $d>0$ and set $\widetilde{d}_N=d_N N^d$.
There exist constants $\epsilon,\delta>0,\alpha>0$, $C\ge1$ independent
of $N$ such that
\[
|L_s^nv|_1\le C\|v\|_b\,
N|b|^\alpha e^{-n\delta\widetilde{d}_N^{\,-1}|b|^{-\alpha}},
\]
for all $v\in F_\theta(\Delta')$, $n\ge1$, and $a\in U_b$.
\end{lemma}

\begin{proof}
Define the sequences $T_{s,n}$, $A_{s,n}$, $B_{s,n}$, $E_{s,n}$
as in Section~\ref{sec-renewal}.
Assuming no approximate eigenfunctions, it follows from
Lemma~\ref{lem-RR} that
$\|T_{s,n}\|_b\le C|b|^\alpha e^{-n\delta\widetilde{d}_N^{\,-1}|b|^{-\alpha}}$.

By truncation of $r$, the operators 
$\|A_{s,n}\|$, $\|B_{s,n}\|$, $\|E_{s,n}\|$ vanish for $n>[q\ln N]$.
For $|a|\le \epsilon N^{-1}$, we compute that
$\|A_{s,n}\|\le C e^{-c'n}$,
$\|B_{s,n}\|\le C Nne^{-c'n}$, and
$\|E_{s,n}\|\le C e^{-c'n}$ where $c'=c-\epsilon$.
The result follows.
\end{proof}

Consequently, $\sum_{n\ge1}|L_s^nv|_1\le C\|v\|_b N
\widetilde{d}_N |b|^{2\alpha}$.  Moreover, it is easy to check that
\[
\|v_s\|_\theta \le C\|v\|_\theta N, \qquad |w_s|_\infty \le C|w|_\infty N.
\]
Hence 
\[
|\hat\rho(s)|\le C \|v\|_\theta |w|_\infty N^3\widetilde{d}_N|b|^{2\alpha},
\]
for all $a\in U_b$ completing Step~(f).
Step (g) is proved as in Section~\ref{sec-exp}, and combining Steps~(b)
and~(g) we obtain
\[
|\rho(t)|\le C \|v\|_{m,\theta}|w|_\infty
\{N^{-\beta}+tN^{-(\beta+1)}+
 N^3\widetilde{d}_N e^{-\epsilon N^{-1}t}+ N^3 \widetilde{d}_N^{\,p+1} t^{-p}  \}.
\]
Set $N=[t/(q\ln t)]$.
For $p$, $q$ sufficiently large, 
$\rho(t)=O((\ln t)^{\beta+1}t^{-\beta})$ as required.

\subsection{Logarithmic factors}
\label{sec-log}

Lemma~\ref{lem-Yn} shows that the decay rates on $r$ and $h$ lead to
a decay rate for $H$.   An alternative approach is to make an assumption
on $H$ (via $Y(n)$) from the outset.   In particular, if we assume that
\[
\mu_Y(r>n)=O(\gamma^n), \qquad \mu_Y(Y(n))=O(n^{-(\beta+1)}),
\]
then we obtain typically the estimate $\rho(t)=O(t^{-\beta})$.
(The proof proceeds by truncating so that $r'=\min\{r,[q\ln N]\}$
and $H'=\min\{H,N\}$, with $U_b$ modified so that
$|a|\le\epsilon N(\ln N)^{-1}$ as in Section~\ref{sec-exp}.) 
With the obvious modifications, we can handle general decay rates
for $\mu_Y(Y(n))$.
Presumably this method gives sharp results, but the assumption
on $H$ is more difficult to verify.

One situation where $Y(n)$ decays slower than $\Delta(n)$ is when
the values of $h$ are constant up the tower.
Write $a_n\sim b_n$ to mean $a_n=O(b_n)$ and $b_n=O(a_n)$.
Suppose that $\mu_Y(r=n)\sim e^{-n}$ and let $h=(n^{-1}e^n)^{1/(\beta+2)}$
on all partition elements $\Delta_{j,\ell}$ with $r(j)=n$.
By definition,
\[
\mu_\Delta(h=(n^{-1}e^n)^{1/(\beta+2)})\sim ne^{-n}.
\]
It follows that $\mu_\Delta(h=n)\sim n^{-(\beta+2)}$ and hence
$\mu_\Delta(\Delta(n))\sim n^{-(\beta+1)}$.
On the other hand, $H=n(n^{-1}e^n)^{1/(\beta+2)}$ on all $Y_j$ with $r(j)=n$
so that a similar calculation gives $\mu_Y(Y(n))\sim(\ln n)^{\beta+1} n^{-(\beta+1)}$ which is one factor of $\ln n$ short of the upper bound
in Lemma~\ref{lem-Yn}.
In this situation, we cannot hope to improve Theorem~\ref{thm-unbounded}.

On the other hand, if we modify the previous example so that $h=e^{n/(\beta+2)}$
on partition elements $\Delta_{j,0}$ with $r(j)=n$ and $h$
is uniformly bounded on the remainder of the tower, then again
$\mu_\Delta(\Delta(n))\sim n^{-(\beta+1)}$, but this time
$\mu_Y(Y(n))\sim n^{-(\beta+1)}$ and typically $\rho(t)=O(t^{-\beta})$.   

Alternatively, suppose that there is a constant $m\ge1$ such that 
for each $j$ with $r(j)\le q\ln N$
 there are at most $m$ values of $\ell<r(j)$ such that
$\|1_{\Delta_{j,\ell}}h\|_\theta\ge N\ln N^{-1}$. 
Then $\|1_{Y_j}H'\|_\theta\le m N+(q\ln N)(N\ln N^{-1})=(m+q)N$.
In this situation, the truncations of $r$ and $h$ automatically achieve the
required truncation of $H$ and so typically $\rho(t)=O(t^{-\beta})$.

\subsection{Flows with unbounded roof function}

Here, we mimic Section~\ref{sec-hyp} but in the context of
Subsection~\ref{sec-semiunb}, taking account of the fact that
$|h'|_\infty=N$.   It is easily verified in~\cite[Lemma~5.4]{Mapp}
that $\bar h$ and $\chi$ inherit a single factor of $N$ from $h'$
and this is compensated for by the fact that $|a|\le\epsilon N^{-1}$.
Proceeding as in Section~\ref{sec-hyp}, we break
the $n$'th term of the series for $\hat\rho(s)$ into 
$I_1+I_2+I_3$ where
\[
I_1,I_2 = O(|b|\gamma_1^k  e^{n\epsilon N\min\{N^{-1},\widetilde{d}_N^{\,-1}|b|^{-\alpha}\}}),
\quad
I_3 = O(|b|^{\alpha+1}N^3\lambda^k e^{-n\delta \widetilde{d}_N^{\,-1}|b|^{-\alpha}}),
\]
where $\lambda>1$.   

To progress further, we modify the definition of $U_b$ so that 
$|a|\le \epsilon N^{-1}\widetilde{d}_N^{\,-1}|b|^{-\alpha}$.  Then
\[
I_1,I_2 = O(|b|^{\alpha+1}N^3\gamma_1^k  e^{n'\epsilon}),
\quad
I_3 = O(|b|^{\alpha+1}N^3\lambda^k e^{-n'\delta}),
\]
where $n'=n\widetilde{d}_N^{\,-1}|b|^{-\alpha}$.
Summing over $n$, we obtain
\[
|\hat\rho(s)|\le C\|v\|_\eta \|w\|_\eta N^3\widetilde{d}_N|b|^{2\alpha+1},
\]
for all $s=a+ib$ with 
$|a|\le \epsilon N^{-1}\widetilde{d}_N^{\,-1}|b|^{-\alpha}$.
Taking $N=t^{1-\epsilon'}$ and $p$ sufficiently large,
we obtain the result claimed in Remark~\ref{rmk-eps}.

\paragraph{Acknowledgements}
This research was supported in part by EPSRC Grant EP/D055520/1
and by a Leverhulme Research Fellowship.
I am grateful to Mark Holland for helpful discussions, and to the 
University of Houston for the use of
e-mail given that pine is inadequately supported on the University
of Surrey network.


\begin{thebibliography}{10}

\bibitem{Aaronson}
J.~Aaronson. \emph{{An Introduction to Infinite Ergodic Theory}}. Math. Surveys
  and Monographs \textbf{50}, Amer. Math. Soc., 1997.

\bibitem{BHM05}
H.~Bruin, M.~Holland and I.~Melbourne. {Subexponential decay of correlations
  for compact group extensions of nonuniformly expanding systems}.
  \emph{Ergodic Theory Dynam. Systems} \textbf{25} (2005) 1719--1738.

\bibitem{Chernov99}
N.~Chernov. Decay of correlations and dispersing billiards. \emph{J. Statist.
  Phys.} \textbf{94} (1999) 513--556.

\bibitem{Chernovsub}
N.~Chernov. {A stretched exponential bound on time correlations for billiard
  flows}.  Preprint, 2006.

\bibitem{ChernovDolgopyatapp}
N.~I. Chernov and D.~Dolgopyat. {Hyperbolic billiards and statistical physics}.
  {Proceedings of International Congress of Mathematicians (Madrid, Spain,
  2006)}.

\bibitem{ChernovYoung00}
N.~Chernov and L.~S. Young. Decay of correlations for {L}orentz gases and hard
  balls. \emph{Hard ball systems and the Lorentz gas}. Encyclopaedia Math. Sci.
  \textbf{101}, Springer, Berlin, 2000, pp.~89--120.

\bibitem{ChernovZhang05}
N.~I. Chernov and H.-K. Zhang. {Billiards with polynomial mixing rates}.
  \emph{Nonlinearity} \textbf{18} (2005) 1527--1553.

\bibitem{ChernovZhang05b}
N.~Chernov and H.-K. Zhang. A family of chaotic billiards with variable mixing
  rates. \emph{Stoch. Dyn.} \textbf{5} (2005) 535--553.

\bibitem{Dolgopyat98a}
D.~Dolgopyat. {On the decay of correlations in Anosov flows}. \emph{Ann. of
  Math.} \textbf{147} (1998) 357--390.

\bibitem{Dolgopyat98b}
D.~Dolgopyat. {Prevalence of rapid mixing in hyperbolic flows}. \emph{Ergodic
  Theory Dynam. Systems} \textbf{18} (1998) 1097--1114.

\bibitem{FMTapp}
M.~J. Field, I.~Melbourne and A.~T{\" o}r{\" o}k. {Stability of mixing and
  rapid mixing for hyperbolic flows}. \emph{Ann. of Math.}  To appear.

\bibitem{FriedmanMartin88}
B.~Friedman and R.~F. Martin. Behavior of the velocity autocorrelation function
  for the periodic {L}orentz gas. \emph{Phys. D} \textbf{30} (1988) 219--227.

\bibitem{Gouezel04a}
S.~Gou{\"e}zel. {Sharp polynomial estimates for the decay of correlations}.
  \emph{Israel J. Math.} \textbf{139} (2004) 29--65.

\bibitem{Hu04}
H.~Hu. Decay of correlations for piecewise smooth maps with indifferent fixed
  points. \emph{Ergodic Theory Dynam. Systems} \textbf{24} (2004) 495--524.

\bibitem{Liverani04}
C.~Liverani. {On contact Anosov flows}. \emph{Ann. of Math.} \textbf{159}
  (2004) 1275--1312.

\bibitem{LiveraniSaussolVaienti99}
C.~Liverani, B.~Saussol and S.~Vaienti. {A probabilistic approach to
  intermittency}. \emph{Ergodic Theory Dynam. Systems} \textbf{19} (1999)
  671--685.

\bibitem{Markarian04}
R.~Markarian. Billiards with polynomial decay of correlations. \emph{Ergodic
  Theory Dynam. Systems} \textbf{24} (2004) 177--197.

\bibitem{Mapp}
I.~Melbourne. {Rapid decay of correlations for nonuniformly hyperbolic flows}.
  \emph{Trans. Amer. Math. Soc.}  To appear.

\bibitem{Pollicott85}
M.~Pollicott. {On the rate of mixing of Axiom A flows}. \emph{Invent. Math.}
  \textbf{81} (1985) 413--426.

\bibitem{Pollicott99}
M.~Pollicott. {On the rate of mixing of Axiom A attracting flows and a
  conjecture of Ruelle}. \emph{Ergod. Th. \& Dynam. Sys.} \textbf{19} (1999)
  535--548.

\bibitem{Ruelle83}
D.~Ruelle. {Flows which do not exponentially mix}. \emph{C. R. Acad. Sci.
  Paris} \textbf{296} (1983) 191--194.

\bibitem{Sarig02}
O.~M. Sarig. {Subexponential decay of correlations}. \emph{Invent. Math.}
  \textbf{150} (2002) 629--653.

\bibitem{Sinai70}
Y.~G. Sina{\u\i}. Dynamical systems with elastic reflections. {E}rgodic
  properties of dispersing billiards. \emph{Uspehi Mat. Nauk} \textbf{25}
  (1970) 141--192.

\bibitem{Young98}
L.-S. Young. Statistical properties of dynamical systems with some
  hyperbolicity. \emph{Ann. of Math.} \textbf{147} (1998) 585--650.

\bibitem{Young99}
L.-S. Young. Recurrence times and rates of mixing. \emph{Israel J. Math.}
  \textbf{110} (1999) 153--188.

\end{thebibliography}
\end{document}